\documentclass[11pt]{article}
\usepackage{amsmath}
\usepackage{amsfonts}
\usepackage{amssymb}
\usepackage{graphicx}
\usepackage{multicol}
\usepackage[usenames]{color}

\usepackage[applemac] {inputenc}
\usepackage[T1]{fontenc}

\usepackage{fullpage}
\def\RR{\rm \hbox{I\kern-.2em\hbox{R}}}
\def\NN{\rm \hbox{I\kern-.2em\hbox{N}}}
\def\ZZ{\rm {{\rm Z}\kern-.28em{\rm Z}}}
\def\CC{\rm \hbox{C\kern -.5em {\raise .32ex \hbox{$\scriptscriptstyle
|$}}\kern
-.22em{\raise .6ex \hbox{$\scriptscriptstyle |$}}\kern .4em}}

\def\nl{\newline}

\def\cM{{\cal M}}
\def\cV{{\cal V}}

\def\cF{{\cal F}}

\def\cD{{\cal D}}

\def\cO{{\cal O}}
\def\R{\mathbb{R}}
\def\N{\mathbb{N}}

\def\Chi{\raise .3ex
\hbox{\large $\chi$}} 
\def\lsima{\hbox{\kern -.6em\raisebox{-1ex}{$~\stackrel{\textstyle<}{\sim}~$}}\kern -.4em}
\def\lsim{\hbox{\kern -.2em\raisebox{-1ex}{$~\stackrel{\textstyle<}{\sim}~$}}\kern -.2em}
\def\gsim{\hbox{\kern -.2em\raisebox{-1ex}{$~\stackrel{\textstyle>}{\sim}~$}}\kern -.2em}

\newcommand{\be}{\begin{equation}}
\newcommand{\ee}{\end{equation}}
\newcommand{\bea}{$$ \begin{array}{lll}}
\newcommand{\eea}{\end{array} $$}
\newcommand{\bi}{\begin{itemize}}
\newcommand{\ei}{\end{itemize}}
\newcommand{\iref}[1]{(\ref{#1})}
\newtheorem{theorem}{Theorem}[section]
\newtheorem{lemma}{Lemma}[section]

\newtheorem{corollary}{Corollary}[section]

\newtheorem{remark}{Remark}[section]

\def\E{\mathbb E}

\def\N{\mathbb N}

\DeclareMathOperator{\essinf}{ess\,inf}

\begin{document}

\title
{Sparse polynomial approximation of
parametric elliptic PDEs \\
Part I: affine coefficients
\thanks{%
Research supported by the European Research Council under grant ERC AdG BREAD.
}
}
\author{ 
Markus Bachmayr, Albert Cohen and Giovanni Migliorati
}

\maketitle
\date{}
\begin{abstract}
We consider the linear elliptic equation $-{\rm div}(a\nabla u)=f$ on some bounded domain $D$, where 
$a$ has the affine form $a=a(y)=\bar a+\sum_{j\geq 1} y_j\psi_j$
for some parameter vector $y=(y_j)_{j\geq 1}\in U=[-1,1]^\N$. We study the summability properties of
polynomial expansions of the solution map $y\mapsto u(y) \in V:= H^1_0(D)$.
We consider both Taylor series 
and Legendre series. Previous results \cite{CDS} 
show that, under a uniform ellipticity assuption, for any $0<p<1$, the 
$\ell^p$ summability of the $(\|\psi_j\|_{L^\infty})_{j\geq 1}$
implies the $\ell^p$ summability of the $V$-norms of the Taylor or Legendre coefficients.
Such results ensure convergence rates $n^{-s}$ of polynomial approximations
obtained by best $n$-term truncation of such series, 
with $s=\frac 1 p-1$ in $L^\infty(U,V)$ or $s=\frac 1 p-\frac 1 2$ in $L^2(U,V)$. In this paper
we considerably improve these results by providing sufficient conditions of
$\ell^p$ summability of the coefficient $V$-norm sequences expressed in terms of the
pointwise summability properties of the $(|\psi_j|)_{j\geq 1}$. The approach in the present paper
strongly differs from that of \cite{CDS}, which is based on individual estimates of the coefficient norms
obtained by the Cauchy formula applied to a holomorphic extension of the solution map.
Here, we use weighted summability estimates, obtained by real-variable arguments.
While the obtained results imply those of \cite{CDS} as a particular case, they lead to a refined
analysis which takes into account the amount of overlap between the supports
of the $\psi_j$. For instance, in the case of disjoint supports, these results
imply that for all $0<p<2$, the $\ell^p$ summability of the coefficient $V$-norm sequences
follows from  the weaker assumption that $(\|\psi_j\|_{L^\infty})_{j\geq 1}$ is $\ell^q$ summable for
$q=q(p):=\frac {2p}{2-p}>p$. We provide a simple analytic example showing that this result is 
in general optimal and illustrate our findings by numerical experiments. 
The analysis in the present paper applies to other types of linear PDEs with 
similar affine parametrization of the coefficients, {and to more general Jacobi polynomial expansions.}

\end{abstract}

\vskip .1in

\section{Introduction}

\subsection{Elliptic PDEs with affine parameter dependence}

The numerical treatment of parametric and stochastic partial differential equations
was initiated in the 1990s, see \cite{GS1,GS2,KL,X} for general references. It
has recently drawn much attention in the case where the number of involved parameters
is very large \cite{BNTT1,BNTT2}, or {\it countably infinite} \cite{CD,CDS,CCS}. 

In this paper, we are interested in the approximation of the elliptic parametric PDE of the form
\be
-{\rm div}(a\nabla u)=f,
\ee
set on a bounded Lipschitz domain $D\subset \R^d$  (where in typical applications, 
$d=1,2,3$), with homogeneous Dirichlet boundary conditions. We consider coefficients having the affine form 
\be
a=a(y)=\bar a+\sum_{j\geq 1}y_j\psi_j,
\label{affine}
\ee
where 
\be
y\in U=[-1,1]^\N.
\ee
The functions $\bar a$ and $\psi_j$ are assumed to be in $L^\infty(D)$.

Well-posedness of this problem in $V=H^1_0(D)$ 
is ensured
for all $y\in U$ by the so-called \emph{uniform
ellipticity assumption}, further referred to as (UEA),
\be
\sum_{j\geq 1}|\psi_j| \leq \bar a-r. \tag{UEA}
\ee
for some $r>0$, which is equivalent to $r\leq a(y)\leq R$ for all $y\in U$ for some $0<r\leq R<\infty$.
This assumption ensures the uniform boundedness of the solution map $y\mapsto u(y)$ from $U$ to $V$.

We are interested in polynomial approximations of the solution map. Such approximations have the
form
\be
u_\Lambda(y):=\sum_{\nu\in \Lambda} v_\nu y^\nu,
\ee
where $\Lambda\subset \cF$ is a finite set of (multi-)indices $\nu=(\nu_j)_{j\geq 1} \in\cF$
and  $y^\nu=\prod_{j\geq 1} y_j^{\nu_j}$.
In the infinite dimensional setting, the index set $\cF$ 
denotes the (countable) set of 
all sequences of nonnegative integers which are 
{\em finitely supported} (i.e. those sequence for which
only finitely many terms are nonzero). Note that the polynomial
coefficients $v_\nu$ are functions in $V$. Such approximations
can be obtained by truncation of infinite polynomial expansions.

Two relevant such expansions are:
\begin{enumerate}
\item The Taylor series
\be
u = \sum_{\nu\in\cF} t_\nu y^\nu, \quad t_\nu:=\frac 1 {\nu !} \partial^\nu u(0),
\ee
where we use the customary multi-index notations 
\be
y^\nu = \prod_{j\geq 1} y_j^{\nu_j} \quad {\rm and}  \quad \nu! = \prod_{j\geq 1} \nu_j!,
\ee
with the convention that $0!=1$.
\item
The Legendre series
\be\label{legendreseries}
u = \sum_{\nu\in\cF} u_\nu L_\nu, \quad u_\nu=\langle u,L_\nu\rangle , \quad L_\nu(y)=\prod_{j\geq 1}L_{\nu_j}(y_j),
\ee
where the univariate Legendre polynomials $L_k$ are normalized in $L^2([-1,1],\frac {dt}2)$, and $\langle \cdot,\cdot\rangle $
denotes the inner product in $L^2(U,V,\mu)$, where $\mu=\bigotimes_{j\geq 1} \frac{dy_j}2$ is the uniform probability measure.
\end{enumerate}

Given such expansions, one natural way of constructing a polynomial approximation
is by best $n$-term truncation, that is, setting
\be
u_n^T(y):=\sum_{\nu\in \Lambda_n^T}t_\nu y^\nu,
\ee
or 
\be
u_n^L(y):=\sum_{\nu\in \Lambda_n^L}u_\nu L_\nu(y),
\ee
where $\Lambda_n^T$ and $\Lambda_n^L$ are the index sets corresponding to the
$n$ largest $\|t_\nu\|_V$ or $\|u_\nu\|_V$, respectively. The convergence rates of
such approximations are governed by the $\ell^p$ summability properties of the sequences 
$(\|t_\nu\|_V)_{\nu\in\cF}$ for $p<1$ and $(\|u_\nu\|_V)_{\nu\in\cF}$ for $p<2$. Indeed, by standard 
application of Stechkin's lemma \cite{De}, such summability properties imply error estimates 
of the form
\be
\|u-u_n^T\|_{L^\infty(U,V)} \leq C(n+1)^{-s}, \quad s:=\frac 1 p-1, \quad C:=\|(\|t_\nu\|_V)_{\nu\in\cF}\|_{\ell^p},
\label{nT}
\ee
and
\be
\|u-u_n^L\|_{L^2(U,V,\mu)} \leq C(n+1)^{-r}, \quad r:=\frac 1 p-\frac 1 2, \quad C:=\|(\|u_\nu\|_V)_{\nu\in\cF}\|_{\ell^p},
\label{nL}
\ee
see \cite{CDS}. { Let us mention that the case of finitely many variables, that is,
when $\psi_j=0$ for $j>d$, typically leads to exponential convergence rates of the form $\exp(-c n^{1/d})$.
Such rates have been obtained in \cite{BNTT1,BNTT2} by application of the above Stechkin estimates
for all $0<p<1$ and tuning the value of $p$ with that of $n$. They can also be directly obtained 
through the available upper bounds on the coefficients $\|t_\nu\|_V$ or $\|u_\nu\|_V$,
as shown in \cite{TWZ}, which brings certain improvements over the previous approach. 
In the present paper we focus our attention on the infinite-dimensional case.}

The above approximation estimates have general implications on the
potential performance of other model reduction techniques. Indeed
we observe that $u_n^T(y)$ and $u_n^L(y)$ belong to fixed $n$-dimensional
subspaces of $V$, spanned by the $(t_\nu)_{\nu\in \Lambda_n^T}$ 
or $(u_\nu)_{\nu\in \Lambda_n^L}$. It thus follows from \iref{nT} that
the best $n$-dimensional model reduction error in the uniform sense, given by
the Kolmogorov $n$-width 
\be
d_n(\cM):=\inf_{\dim(V_n)=n}\sup_{v\in \cM}\min_{w\in V_n}\|v-w\|_V, \quad \cM:=u(U)=\{u(y)\; : \; y\in U\},
\ee
decays at least as fast as $n^{-s}$. Likewise, it follows from \iref{nL} that
if $y$ is uniformly distributed in $U$, the best $n$-dimensional model reduction error in the mean-square sense
\be
\inf_{\dim(V_n)=n}\E(\min_{w\in V_n}\|u(y)-w\|_V^2),
\ee
attained when $V_n$ is the span of the $n$ first $V$-principal components,
decays at least as fast as $n^{-2r}$. The above estimates govern the convergence
rate of reduced basis (RB) methods and proper orthogonal decomposition (POD) methods, respectively,
see \cite{CD}.

\subsection{Existing results}

Summability results for the sequences 
$(\|t_\nu\|_V)_{\nu\in\cF}$ and $(\|u_\nu\|_V)_{\nu\in\cF}$ have been established in \cite{CDS} under the so-called \emph{uniform
ellipticity assumption}
\be\label{uea}
\sum_{j\geq 1}|\psi_j| \leq \bar a-r. \tag{UEA}
\ee
for some $r>0$, which is equivalent to $r\leq a(y)\leq R$ for all $y\in U$ for some $0<r\leq R<\infty$.
This assumption ensures the uniform boundedness of the solution map $y\mapsto u(y)$ from $U$ to $V$.
These results can be summarized as follows.

\begin{theorem}
\label{firsttheo}
Assume that (UEA) holds. For any $0<p<1$, 
if $(\|\psi_j\|_{L^\infty})_{j\geq 1} \in \ell^p(\N)$ then $(\|t_\nu\|_V)_{\nu\in\cF}\in \ell^p(\cF)$ and 
$(\|u_\nu\|_V)_{\nu\in\cF}\in \ell^p(\cF)$.
\end{theorem}

These results have been extended to a large range of linear or nonlinear
parametric PDEs \cite{CD,CCS} where $y$ is again ranging in the infinite dimensional box $U$.
They strongly rely on the holomorphy of the solution map $y\mapsto u(y)$ in each
variable $y_j$. While they yield provable algebraic approximation rates for parametric PDEs 
in the infinite dimensional framework,
these results are not fully satisfactory for several reasons.

First, {they are confined to the case $0<p<1$. Therefore, no summability
property of the Taylor or Legendre coefficients can be deduced from these results in
cases where $(\|\psi_j\|_{L^\infty})_{j\geq 1}$ does not have this type of $\ell^p$
summability. Note also that, in the case of Legendre series, 
convergence rates of best $n$-term approximation can be derived as soon as 
$(\|u_\nu\|_V)_{\nu\in\cF}\in \ell^p(\cF)$ for some $p<2$, in view of \iref{nL}.
A legitimate objective is therefore to identify conditions on $(\psi_j)_{j\geq 1}$ that
govern the $\ell^p$ summability of $(\|u_\nu\|_V)_{\nu\in\cF}$
in the regime $1\leq p\leq 2$.}

Second and more important, the summability conditions imposed on the $\|\psi_j\|_{L^\infty}$ in this result becomes
quite strong and artificial in the case where the supports of these functions do not overlap too much.
As a relevant example, consider the case where the $(\psi_j)_{j\geq 1}$ are a wavelet basis on the domain $D$. In this case
it is more natural to denote the elements of such bases by $\psi_\lambda$, where $\lambda$ concatenates
the space and scale indices, following
the usual terminology such as in \cite{Co}, and using the the notation $l=|\lambda|$ for the scale level.
At a given scale level $l\geq 0$, there
are $\cO(2^{dl})$ wavelets and their supports have finite overlap in the sense that, for all $x\in D$,
\be
\#\{\lambda \,  : \, |\lambda|=l\quad {\rm and}\quad \psi_\lambda(x)\neq 0\} \leq M,
\label{finiteoverlap}
\ee
for some fixed $M$ independent of $l$, {where $\# S$ denotes the cardinality of the set $S$}. It is well known that the geometric rate of decay of the wavelet contributions as the scale level grows reflects 
the amount of smoothness in the expansion.
It is thus natural to study the situation where
\be
\|\psi_\lambda\|_{L^\infty}= c_l:= C2^{-\alpha l},\quad |\lambda|=l,
\ee
for some given $\alpha>0$. In the affine case, assuming that $\bar a$ and the $\psi_\lambda$ are sufficiently smooth, this means that the diffusion coefficient $a$ is uniformly bounded in the H\"older space $C^\alpha(D)$ independently of $y$. In this case, the uniform ellipticity assumption is ensured 
provided that
\be
M\sum_{l\geq 0} c_l \leq \bar a_{\min} -r, \quad \text{where}\quad \bar a_{\min} := \essinf_D \bar a >0,
\ee
which holds if $C\leq M^{-1}(\bar a_{\min} -r)(1-2^{-\alpha})$ for the above constant $C$. Note that we can take $\alpha>0$
arbitrarily small, up to taking $C$ small enough. The uniform ellipticity assumption ensures that the solution map
belongs to $L^\infty(U,V)$ and therefore to $L^2(U,V,\mu)$, and thus that the Legendre coefficient sequence
$(\|u_\nu\|_V)_{\nu\in \cF}$ belongs to $\ell^2(\cF)$. However, if we want to use the above mentioned results to prove
$\ell^p$ summability of this sequence for smaller values of $p$, we are enforced to assume that 
the $\|\psi_\lambda\|_{L^\infty}$ are summable over all indices, which equivalently means that
\be
\sum_{l\geq 0} 2^{dl} 2^{-\alpha l}<\infty,
\ee
and therefore that $\alpha >d$. This is highly unsatisfactory since it excludes diffusion coefficients with low
order of smoothness, especially when the spatial dimension $d$ is large.

\subsection{Main results and outline}

The above example reveals a gap in the currently available analysis: $\ell^2$ summability of $(\|v_\nu\|_V)_{\nu\in \cF}$ can be obtained 
under mild assumptions, while proving $\ell^p$ summability for $p<2$
by the existing results immediately imposes much stronger assumptions (in the sense of the required decay of $\|\psi_j\|_{L^\infty}$ as $j\to \infty$).
In this paper, we propose a new analysis which allows us to remove this gap.

The main results of the present paper can be summarized as follows.

\begin{theorem}
\label{maintheo}
Assume that \iref{uea} holds and that for $0<p<2$ and $q=q(p):=\frac {2p}{2-p}$,
there exists a sequence $\rho=(\rho_j)_{j\geq 1}$ with $\rho_j>1$, such that
\be
\sum_{j\geq 1} \rho_j|\psi_j| \leq \bar a-s,
\ee
for some $s>0$, and such that $\bigl(\rho_j^{-1} \bigr)_{j\geq 1}$ belongs to $\ell^q(\N)$.
Then the sequences
$\bigl(\|t_\nu\|_V \bigr)_{\nu\in \cF}$ and 
$\bigl(\|u_\nu\|_V\bigr)_{\nu\in \cF}$ belong to $\ell^p(\cF)$.
\end{theorem}

This theorem constitutes a particularly strong improvement over Theorem \ref{firsttheo} when 
the supports of the functions $(\psi_j)_{j\geq 1}$ do not overlap much. For example, in the case
where these supports are disjoint we can take weights such that
$\rho_j^{-1}\sim \|\psi_j\|_{L^\infty}$ as $\|\psi_j\|_{L^\infty}$ tends to $0$. Therefore,
in this case, for all $0<p<2$ the $\ell^p$ summability of 
$(\|t_\nu\|_V)_{\nu\in \cF}$ and $(\|u_\nu\|_V)_{\nu\in \cF}$
follows from the assumption that $(\|\psi_j\|_{L^\infty})_{j\geq 1}$ is $\ell^q$ summable for
$q=q(p):=\frac {2p}{2-p}$. Note that $q(p)>p$ and that
\be
\lim_{p\to 2} q(p)= +\infty,
\ee
which shows that almost no decay of $(\|\psi_j\|_{L^\infty})_{j\geq 1}$
is required as $p$ gets closer to $2$. Similar improvements can be obtained for
other types of families $(\psi_j)_{j\geq 1}$, such as wavelets.

Let us also mention that,
while we focus on the diffusion equation, inspection of proofs
reveals that the main results can be extended to other types of linear elliptic or parabolic PDEs with 
similar affine dependence of the coefficients. { Specific examples are given at the end of \S 2.}

The remainder of this paper is organized as follows. In \S 2, we establish Theorem \ref{maintheo}
in the case of Taylor coefficients. The approach strongly differs from that used in \cite{CDS}
for proving Theorem \ref{firsttheo} which is based on establishing individual estimates of the $\|t_\nu\|_V$
by Cauchy's formula applied to to a holomorphic extension of the solution map. Instead,
we use weighted summability estimates, obtained by real-variable arguments. It is worth mentioning
that these estimates also imply the individual estimates from \cite{CDS}.

In \S 3, we establish Theorem \ref{maintheo} in the case of Legendre coefficients. The approach
is again based on weighted summability estimates, obtained by combining the ideas
developed for the Taylor coefficients with Rodrigues' formula. {Furthermore, we give an extension of our results to Jacobi polynomial expansions by an analogous argument.}

We give in \S 4 several examples of applications, corresponding to different types
of support properties for the $(\psi_j)_{j\geq 1}$ and we discuss in each case the improvements
over Theorem \ref{firsttheo}. In particular, we show that our results are sharp in
the case of disjoint or finitely overlapping supports, in the sense that $\ell^q$ summability
of $(\|\psi_j\|_{L^\infty})_{j\geq 1}$ with $q:=q(p)$ is generally necessary to achieve $\ell^p$ summability of 
$(\|t_\nu\|_V)_{\nu\in \cF}$. {While we focus on the infinite-dimensional framework, the obtained
results can also be used in the case of finitely many variables. In particular, one may use them to obtain
improved exponential convergence rates, by following either a Stechkin-based approach as in \cite{BNTT1,BNTT2} or a more direct approach as in \cite{TWZ}.}

Finally, we give in \S 5 the results of several numerical tests 
evaluating the sharpness of the obtained results for the various types of $(\psi_j)_{j\geq 1}$
considered in \S 4. 

\section{Summability of Taylor coefficients}

We start from an analysis of Taylor coefficients, where we use the 
following alternative expression of the uniform ellipticity assumption, further referred to as (UEA*),
\be\label{uea1}
   \Biggl\|  \frac{\sum_{j\geq 1} | \psi_j | }{ \bar a} \Biggr\|_{L^\infty} < 1.\tag{UEA*}
\ee
Indeed, if \eqref{uea} holds with some $r>0$, then we also have
\be
   \Biggl\|  \frac{\sum_{j\geq 1} | \psi_j | }{ \bar a} \Biggr\|_{L^\infty} \leq \biggl\| \frac{\bar a - r}{\bar a} \biggr\|_{L^\infty} \leq 1 - \frac{r}{\| \bar a\|_{L^\infty}} < 1,
\ee
and it is also easily checked that \iref{uea1} implies \iref{uea} for a certain $r$.

Similar arguments as in \cite{CDS} show that under such an assumption, the partial derivatives
$\partial^\nu u(y)$ are well defined for each $\nu\in \cF$ as elements of $V$ for each $y\in U$.
They can be computed by applying the operator $\partial^\nu$ in the $y$ variable to the 
variational formulation
\be
\int_D a(y)\nabla u(y)\cdot\nabla v\, dx=\langle f,v\rangle_{V',V}, \quad v\in V,
\ee
which gives, for all $\nu\neq 0$,
\be
\int_D a(y) \nabla \partial^\nu u(y) \cdot \nabla v \,dx = - \sum_{j\in \operatorname{supp} \nu } \nu_j\int_D \psi_j \nabla \partial^{\nu-e_j} u(y)  \cdot \nabla v\,dx , \quad  v\in V,
\label{partialrecurs}
\ee
where 
\be
e_j=(0,\dots,0,1,0,\dots)=(\delta_{i,j})_{i\geq 1}
\ee
is the Kronecker $\delta$ sequence with $1$ at position $j$. For the Taylor coefficients, this yields
\be
  \int_D \bar a \nabla t_\nu \cdot \nabla v \,dx = - \sum_{j\in \operatorname{supp} \nu } \int_D \psi_j \nabla t_{\nu - e_j} \cdot \nabla v\,dx , \quad  v\in V,
  \label{taylorrecurs}
\ee
when $\nu\neq 0$ and 
\be
\int_D \bar a \nabla t_0 \cdot \nabla v \,dx =\langle f,v\rangle_{V',V}, \quad v\in V.
\ee
We shall make use of the norm defined by
\be
\|v\|_{\bar a}^2:=\int_D \bar a|\nabla v|^2 \, dx.
\ee
This norm is obviously equivalent to the $V$-norm with
\be
\bar a_{\min} \|v\|_V^2 \leq \|v\|_{\bar a}^2 \leq \|\bar a\|_{L^\infty} \|v\|_V^2.
\label{normequiv}
\ee
 The following results shows that under \eqref{uea1}, the energy of Taylor coefficients
decays geometrically with the total order of differentiation. Its proofs uses a technique introduced
in \cite{CCDS}.
Here we use the standard notation
\be
|\nu|:=\sum_{j\geq 1} \nu_j,
\ee
to denote the total order.

\begin{lemma}\label{lmm:contr}
If \eqref{uea1} holds, then for $\sigma:=\frac {\theta}{2-\theta}<1$ with
\be\label{thetadef}
 \theta=\Biggl \| \frac{\sum_{j\geq 1}|\psi_j|}  {\bar a}\Biggr \|_{L^\infty} <1  ,
\ee
we have for all $k\geq 1$,
\be\label{eq:contr}
\sum_{|\nu|=k} \|t_\nu\|^2_{\bar a} \leq   \sigma \sum_{|\nu|=k-1} \|t_\nu\|^2_{\bar a}\,.
\ee
\end{lemma}

\noindent
{\bf Proof:} For $\nu \in \cF$, we define
\be
  d_\nu := \int_D \bar a |\nabla t_\nu |^2 \,dx ,\quad  d_{\nu,j} := \int_D | \psi_j | | \nabla t_\nu |^2\,dx \,.
\ee
Then by \eqref{uea1},
\be\label{eq:dnujestimate}
  \sum_{j \geq 1} d_{\nu,j} \leq  \theta d_\nu .
\ee
Furthermore, since for $\nu\neq 0$ the Taylor coefficients satisfy
\be
  \int_D \bar a \nabla t_\nu \cdot \nabla v \,dx = - \sum_{j\in \operatorname{supp} \nu } \int_D \psi_j \nabla t_{\nu - e_j} \cdot \nabla v\,dx , \quad  v\in V,
\ee
we also have
\be 
   d_\nu \leq \sum_{j\in \operatorname{supp} \nu } \int_D |\psi_j| |\nabla t_{\nu-e_j}| |\nabla t_\nu| \,dx \leq
     \frac12  \sum_{j\in \operatorname{supp} \nu }  \bigl( d_{\nu-e_j,j} + d_{\nu,j} \bigr) ,
\ee
where we have used Young's inequality. Thus, by \eqref{eq:dnujestimate},
\be
  \biggl( 1 - \frac\theta2\biggr) d_\nu \leq \frac12 \sum_{j\in \operatorname{supp} \nu }  d_{\nu-e_j,j} .
\ee
Summing over $|\nu|=k$, we obtain
\be
\biggl( 1 - \frac\theta2\biggr) \sum_{|\nu|=k}d_\nu \leq \frac12 \sum_{|\nu|=k}\sum_{j\in \operatorname{supp} \nu }  d_{\nu-e_j,j}=
 \frac12 \sum_{|\nu|=k-1}\sum_{j\geq 1 }  d_{\nu,j},
\ee
and using again \eqref{eq:dnujestimate}, we arrive at \eqref{eq:contr}.\hfill$\Box$
\nl

Note that as an immediate consequence of \eqref{eq:contr},
\be\label{eq:sumbound0}
  \sum_{\nu\in\cF} \|t_\nu\|^2_{\bar a} \leq    \frac{  \| t_0\|^2_{\bar a}  }{  1 - \sigma  } =\frac{  2 - \theta   }{ 2 - 2\theta }   \| t_0\|^2_{\bar a}< \infty,
\ee
and therefore by \iref{normequiv} and the Lax-Milgram stability estimate for $t_0=u(0)$,
\be\label{eq:sumbound}
  \sum_{\nu\in\cF} \|t_\nu\|^2_{V} \leq    \frac{ ( 2 - \theta ) \|\bar a\|_{L^\infty} }{ (2 - 2\theta) \bar a_{\min} }    \| t_0\|^2_{V}    
  \leq  \frac{ ( 2 - \theta ) \|\bar a\|_{L^\infty} }{ (2 - 2\theta) \bar a_{\min}^3 }    \| f\|^2_{V'}    
  < \infty \,.
\ee
We are now ready to state the main result of this section which is a direct consequence of the above observations.

\begin{theorem}\label{thm:weightedsum}
Assume that for some sequence $\rho=(\rho_j)_{j\geq 1}$ of positive weights, we have
the weighted uniform ellipticity assumption
\be\label{suea}
\delta:= \Biggl\|  \frac{\sum_{j\geq 1} \rho_j |\psi_j|}{\bar a} \Biggr\|_{L^\infty}  < 1.
\ee
Then the sequence $(\rho^\nu\|t_\nu\|_V)_{\nu\in \cF}$ is $\ell^2$ summable, with
\be
\sum_{\nu\in\cF} (\rho^\nu\|t_\nu\|_V)^2 \leq C <\infty,
\ee
where
\be
C=C(\bar a,f,\delta)=\frac{ ( 2 - \delta ) \|\bar a\|_{L^\infty} }{ (2 - 2\delta) \bar a_{\min}^3 }    \| f\|^2_{V'}.
\label{constant}
\ee
\end{theorem}

\noindent
{\bf Proof:} 
We observe that the weighted UEA is equivalent to \eqref{uea1} for the rescaled coefficient
$a_\rho(y):=a(D_\rho y)$, where $D_\rho y:=(\rho_j y_j)_{j\geq 1}$. So we obtain with Lemma \ref{lmm:contr} and \eqref{eq:sumbound} that
\be
\sum_{\nu\in\cF} \|t_{\rho,\nu}\|^2 <\infty,
\ee
where 
\be
t_{\rho, \nu}=\frac 1 {\nu !} \partial^\nu u_\rho(0)= \rho^\nu t_\nu, \quad u_\rho(y)=u(D_\rho y).
\ee
The result follows. \hfill $\Box$
\nl

As a consequence we obtain the following summability result.

\begin{corollary}\label{sumresult}
If for some sequence $\rho=(\rho_j)_{j\geq 1}$ with $\rho_j>1$, $j\in \N$, we have
the weighted uniform ellipticity assumption \eqref{suea} and if the sequence
$\bigl(\rho_j^{-1} \bigr)_{j\geq 1}$ belongs to $\ell^q(\N)$ with $q=q(p):=\frac {2p}{2-p}$ for some $p<2$, then 
the sequence $\bigl(\|t_\nu\|_V \bigr)_{\nu\in \cF}$ is $\ell^p$ summable.
\end{corollary}

\noindent
{\bf Proof:} By H\"older's inequality,
\be
  \sum_{\nu\in\cF} \| t_\nu \|_V^p  \leq  \Bigl( \sum_{\nu\in\cF}  \rho^{2\nu} \| t_\nu \|_V^2  \Bigr)^{p/2}  \Bigl( \sum_{\nu\in\cF} \rho^{-\frac{2p}{2-p} \nu}  \Bigr)^{(2-p)/2}.
\ee
Moreover,
\be
 \sum_{\nu\in\cF} \rho^{- \frac{2p}{2-p} \nu}  = \prod_{j\geq 1} \Bigl( \sum_{k=0}^\infty \rho_j^{-qk}  \Bigr) = \prod_{j\geq 1} ( 1- \rho^{-q}_j)^{-1},
\ee
where the latter product converges precisely when $(\rho^{-1}_j)_{j\geq 1} \in\ell^q$.
The statement thus follows from Theorem \ref{thm:weightedsum}.\hfill$\Box$

\begin{remark}
\label{remcor}
As a trivial consequence of Theorem \ref{thm:weightedsum}, using the fact that the $\ell^2$ norm dominates the $\ell^\infty$ norm, 
we also retrieve the estimate
\be
\|t_\nu\|_V\leq C\rho^{-\nu},
\ee
for any sequence $\rho$ such that \iref{suea} holds, where $C$ is the square root of the constant in \iref{constant}.
This estimate was established in \cite{CDS}, with a different
constant, by a complex variable argument, namely invoking the holomorphy of $y\mapsto u(y)$ on polydiscs of the form
$\otimes_{j\geq 1} \{|z_j|\leq \rho_j\}$. One advantage of this individual estimate is that one may choose to optimize 
over all possible sequences $\rho$, which yields
\be
\|t_\nu\|_V\leq C\inf \rho^{-\nu},
\ee
where the infimum is taken over all sequences $\rho=(\rho_j)_{j\geq 1}$ of numbers larger than $1$ such that \iref{suea} holds.
The proof of Theorem \ref{firsttheo} in \cite{CDS} for the Taylor coefficients is based on using the above estimate, which
amounts to selecting a different sequence $\rho=\rho(\nu)$ for each $\nu\in\cF$.
However, we show in \S 4 that in several relevant cases better results can 
be obtained by using Corollary \ref{sumresult}.
\end{remark}

{
The above analysis exploits the particular structure of the linear diffusion equation
and affine dependence of $a$ in the variables $y_j$, in particular when deriving
the recursion equation \iref{taylorrecurs} which subsequently leads to the $\ell^2$ estimate
of Lemma \ref{lmm:contr}. It can be extended in a natural manner to other types of
equations with similar properties. Here are two simple examples.
\begin{enumerate}
\item
The fourth order equation $\Delta(a\Delta u)=f$ set on a Lipschitz domain $D$ with homogeneous boundary condition $u=\frac {\partial u}{\partial n}=0$ on $\partial D$
and $a$ having the affine form \iref{affine}. In this case the relevant space is 
$V:=H^2_0(D)$, the closure of $\cD(D)$ in $H^2(D)$, and the recursion is
\be
  \int_D \bar a \Delta t_\nu \Delta v \,dx = - \sum_{j\in \operatorname{supp} \nu } \int_D \psi_j \Delta t_{\nu - e_j}  \Delta v\,dx , \quad  v\in V.
  \label{taylorrecursbilap}
\ee
By taking $v=t_\nu$ we obtain Lemma \ref{lmm:contr} for the norm $\|v\|_{\bar a}:=\bigl(\int_D  \bar a |\Delta v|^2\bigr)^{1/2}$, which
is equivalent to $\|v\|_V=\|\Delta v\|_{L^2}$.
\item
The parabolic equation $\partial_t u-{\rm div}(a\nabla u)=f$ set on $]0,T[\times D$ with 
homogeneous boundary condition $u=0$ on $]0,T[\times \partial D$ and initial condition $u_{|t=0}=u_0\in L^2(D)$.
In this case, the relevant solution space is $\cV=L^2(]0,T[;H^1_0(D))\cap H^1(]0,T[;H^{-1}(D))$. 
We refer to
\cite{SS} for the corresponding space-time variational formulation and the analysis of its well-posedness. For $\nu\neq 0$
we find that $t_\nu$ satisfies a parabolic PDE with initial value $t_{\nu|t=0}=0$ and recursive variational form
given by
\be
\int_{0}^T\int_D (\partial_t t_\nu \, v +\bar a \nabla t_\nu \cdot \nabla v)  \,dx\, dt= - \sum_{j\in \operatorname{supp} \nu }
\int_{0}^T \int_D \psi_j \nabla t_{\nu - e_j} \cdot \nabla v\,dx , \quad  v\in V.
\ee
By taking $v=t_\nu$ we obtain Lemma \ref{lmm:contr} for the norm $\|v\|_{\bar a}:=\bigl( \int_0^T\int_D  \bar a |\nabla v|^2\bigr)^{1/2}$.
Note that this norm is not equivalent to the full norm of $\cV$ but only to that of $V=L^2(]0,T[;H^1_0(D))$.
\item
Similar results hold when $a$ is of symmetric tensor type, with $\bar a$ non-degenerate uniformly over $D$. 
We obtain Lemma 2.1 for the norm $\|v\|_{\bar a}:=\bigl(\int_D  \bar a |\nabla v|^2 \bigr)^{1/2}$,
under a uniform ellipticity assumption taking the following form, 
\be
\theta:=\left \| \frac{\sum_{j\geq 1} \|\psi_j\|_2}{\lambda_{\min} (\bar a)}  \right \|_{L^\infty} <1.
\ee
analogous to \iref{thetadef}, with $\|\psi_j\|_2$ and $\lambda_{\min} (\bar a)$ denoting the functions 
$x\mapsto \|\psi_j(x)\|_2$ (with $\|\cdot\|_2$ the spectral norm) and $x\mapsto \lambda_{\min} (\bar a)$. Note that this assumption is equivalent
to \iref{thetadef}, that is, to (UEA*) or (UEA) in the scalar case. Likewise we obtain the subsequent results
under the assumption that the sequence $(\rho_j)_{j\geq 1}$ satisfies
$\left \| \frac{\sum_{j\geq 1} \rho_j\|\psi_j\|_2}{\lambda_{\min} (\bar a)}  \right \|_{L^\infty} <1.$
\end{enumerate}
}

\section{Summability of Legendre coefficients}

In this section, we show that the summability properties of Corollary \ref{sumresult} hold also for the Legendre coefficients of $u$.

\begin{theorem}\label{thm:legendreweightedsum}
If for some sequence $\rho=(\rho_j)_{j\geq 1}$ with $\rho_j \geq 1$, $j\in\N$, we have
the weighted uniform ellipticity assumption
\be\label{suea2}
\delta:=\Biggl\|  \frac{\sum_{j\geq 1} \rho_j |\psi_j| }{\bar a} \Biggr\|_{L^\infty} < 1,
\ee
then with $a_\nu := \prod_{j\geq 1} \sqrt{2\nu_j + 1}$, the sequence $(a_\nu^{-1} \rho^\nu\|u_\nu\|_V)_{\nu\in \cF}$ is $\ell^2$ summable, that is,
\be
\label{weightedell2}
   \sum_{\nu \in \cF} (a_\nu^{-1} \rho^\nu\|u_\nu\|_V)^2 \leq C < \infty,
\ee
where
\be
 C = C(\bar a, f, \delta) := \frac{(2-\delta)(1+\delta) \|\bar a\|^2_{L^\infty}  \| f\|_{V'}^2  }{ 2(1-\delta)^4 \bar{a}_{\min}^4 }  .
\ee
\end{theorem}

\noindent
{\bf Proof:} For $y, z \in U$, we set $T_y z := \bigl( y_j +  ( 1- |y_j|) \rho_j z_j \bigr)_{j\geq 1}$. Then for $w_y(z) := u(T_y z)$, we have
\be\label{eq:tayloridentity}
  \partial^\nu w_y(0) =  \Bigl( \prod_{j \geq 1}  (1 - |y_j|)^{\nu_j}  \Bigr)  \rho^\nu \partial^\nu u(y) \,.
\ee
Let us fix $y\in U$ and set  $\bar a_y :=a(y)= \bar a + \sum_{j\geq 1} y_j \psi_j$ and $\psi_{y,j} := (1 - |y_j|) \rho_j \psi_j$. Then $w_y$ is the solution of
\be\label{pdewy}
   -\operatorname{div} \Bigl[ \Bigl(  \bar a_y  + \sum_{j\geq 1} z_j \psi_{y,j} \Bigr) \nabla w_y(z) \Bigr] = f \,.
\ee
Applying Lemma \ref{lmm:contr} with the modified $y$-dependent coefficients in \eqref{pdewy}, for the Taylor coefficients $ t_{y,\nu} := (\nu!)^{-1} \partial^\nu w_y(0)$ of $w_y$ we obtain
\be
  \sum_{|\nu|=k} \|  t_{y,\nu}\|^2_{\bar a_y} \leq  \sigma_y  \sum_{|\nu|=k-1} \| t_{y,\nu}\|^2_{\bar a_y},\quad  \sigma_y = \frac{\theta_y}{2-\theta_y},\quad
   \theta_y = \Biggl\| \frac{\sum_{j\geq 1} | \psi_{y,j}|}{\bar a_y} \Biggr\|_{L^\infty} \,.
\ee
Since $\rho_j \geq 1$,
\be
  \theta_y \leq  \Biggl\|  \frac{\sum_{j\geq 1} \rho_j | \psi_{j}| - \sum_{j\geq 1} \rho_j |y_j| |\psi_j| }{\bar a - \sum_{j\geq 1} |y_j| |\psi_j| }\Biggr\|_{L^\infty} 
    \leq  \Biggl\|  \frac{\sum_{j\geq 1} \rho_j | \psi_{j}|  }{\bar a } \Biggr\|_{L^\infty}   = \delta < 1,
\ee
and thus $\sigma_y \leq  \delta/(2-\delta)<1$.
Consequently, as in \eqref{eq:sumbound},
\be\label{eq:legendresumbound}
 \sum_{\nu\in\cF} \|  t_{y,\nu}\|^2_{\bar a}  \leq   \| \bar a_y^{-1} \bar a  \|_{L^\infty}  \sum_{\nu\in\cF} \| t_{y,\nu}\|^2_{\bar a_y}  \leq 
     \| \bar a_y^{-1} \bar a  \|_{L^\infty}  \frac{2-\delta }{ 2 - 2\delta }  \|  t_{y,0} \|^2_{ \bar{a}_y }  \,,
\ee
where 
\be
  \|  t_{y,0} \|^2_{ \bar{a}_y }   \leq  \| \bar{a}_y \|_{L^\infty}   \| \bar{a}_y^{-2}\|_{L^\infty}  \| f\|_{V'}^2 .
\ee
We also have 
\be
     \| \bar a_y^{-1} \bar a  \|_{L^\infty}    \| \bar{a}_y \|_{L^\infty}   \| \bar{a}^{-2}_y \|_{L^\infty}  \leq  \| \bar a\|_{L^\infty} \Bigl\| \bar{a} + \sum_{j\geq 1} | \psi_j|  \Bigr\|_{L^\infty}  \Bigl\|  \Bigl( \bar{a} - \sum_{j\geq 1} | \psi_j|  \Bigr)^{-3} \Bigr\|_{L^\infty}  \leq \frac{(1+\delta) \|\bar a\|^2_{L^\infty}}{ (1-\delta)^3 \bar{a}_{\min}^3 }
\ee
by \eqref{suea2}. Altogether, we obtain
\be
  \sum_{\nu\in\cF} \|t_{y,\nu}\|^2_V  \leq \bar{a}_{\min}^{-1}  \sum_{\nu\in\cF} \|  t_{y,\nu}\|^2_{\bar a} \leq \frac{(2-\delta)(1+\delta) \|\bar a\|^2_{L^\infty}  \| f\|_{V'}^2  }{ (2 -2\delta) (1-\delta)^3 \bar{a}_{\min}^4 }  = C < \infty.
\ee 

With the present normalization as in \eqref{legendreseries}, the Legendre polynomials satisfy the Rodrigues' formula
\be
  L_\nu(y) = \prod_{j\geq 1} \partial^{\nu_j}_{y_j}\biggl(  \frac{\sqrt{2 \nu_j + 1}}{\nu_j! \,2^{\nu_j}} ( y_j^2 - 1)^{\nu_j} \biggr) .
\ee
As a consequence, for the Legendre coefficients of $u$ we obtain
\begin{equation}\label{rodriguesapplication}
\begin{aligned}
   u_\nu &= \int_U u(y) L_\nu(y)\,d\mu(y) = \int_U u(y)  \prod_{j\geq 1} \partial^{\nu_j}_{y_j}\biggl(  \frac{\sqrt{2 \nu_j + 1}}{\nu_j! \,2^{\nu_j}} ( y_j^2 - 1)^{\nu_j} \biggr) \,d\mu(y)   \\ 
     &= \Bigl( \prod_{j \geq 1} \sqrt{2 \nu_j + 1} \Bigr) \int_U \frac{1}{\nu!}\partial^\nu u(y) \, \prod_{j\geq 1} \frac{ ( 1 - y_j^2)^{\nu_j} }{2^{\nu_j}} \,d\mu(y) .
\end{aligned}
\end{equation}
Hence, by \eqref{eq:tayloridentity} and \eqref{eq:legendresumbound},
\begin{equation}\label{legendreestimate}
\begin{aligned}
  \sum_{\nu\in\cF} \Bigr(\prod_{j\geq 1} (2 \nu_j + 1)  \Bigr)^{-1} \rho^{2\nu} \| u_\nu\|_V^2 & \leq \sum_{\nu\in\cF} \rho^{2\nu} \int_U \Bigl\| \frac{1}{\nu!}\partial^\nu u(y) \Bigr\|_V^2 \, \prod_{j\geq 1} \frac{ ( 1 - y_j^2)^{2\nu_j} }{2^{2\nu_j}} \,d\mu(y)   \\
   & = \int_U \sum_{\nu\in\cF} \rho^{2\nu} \Bigl\| \frac{1}{\nu!}\partial^\nu u(y) \Bigr\|_V^2 \prod_{j\geq 1} (1 - |y_j|)^{2\nu_j}\frac{ ( 1 + |y_j|)^{2\nu_j} }{2^{2\nu_j}}   \,d\mu(y)  \\  
   &  \leq  \int_U \sum_{\nu\in\cF}  \Bigl\| \frac{1}{\nu!}\partial^\nu w_y(0) \Bigr\|_V^2 \,d\mu(y) =  \int_U \sum_{\nu\in\cF} \|t_{y,\nu}\|_V^2 \,d\mu(y)   \leq   C ,
\end{aligned}
\end{equation}
which completes the proof.\hfill$\Box$

\begin{corollary}
\label{corleg}
If for some sequence $\rho=(\rho_j)_{j\geq 1}$ with $\rho_j>1$, $j\in \N$, we have
the weighted uniform ellipticity assumption \eqref{suea2}, if the sequence
$\bigl(\rho_j^{-1} \bigr)_{j\geq 1}$ belongs to $\ell^q$ with $q=\frac {2p}{2-p}$ for a $p<2$,
then the sequence $\bigl(\|u_\nu\|_V \bigr)_{\nu\in \cF}$ is $\ell^p$ summable. 
\end{corollary}

\noindent
{\bf Proof:} We obtain the statement from Theorem \ref{thm:legendreweightedsum} by proceeding exactly as in the proof of Corollary \ref{sumresult}.
In this case we need
\be 
  \biggl(\prod_{j\geq1} \sqrt{2\nu_j+1}\, \rho_j^{-\nu_j} \biggr)_{\nu\in\cF} \in \ell^q, 
\ee 
which, since $\rho_j>1$, holds precisely when $(\rho_j^{-1})_{j\geq 1} \in \ell^q$, since
\be
 \sum_{\nu\in\cF} \rho^{- q \nu} \prod_{j\geq1} (2\nu_j+1)^{q/2} = \prod_{j\geq 1} \Bigl( \sum_{k=0}^\infty \rho_j^{-qk}(1+2k)^{q/2}  \Bigr),
\ee
and since, by using the fact that $\|(\rho_j^{-1})_{j\geq 1}\|_{\ell^\infty}<1$, we find that the sum in
each factor of the above product converges and is bounded by $1+C\rho_j^{-q}$ for some fixed $C$.
\hfill$\Box$

\begin{remark}
\label{remleg}
Similar to the case of Taylor coefficients, we can also derive from
Theorem \ref{thm:legendreweightedsum} the individual estimate
\be
\|u_\nu\|_V\leq Ca_\nu \rho^{-\nu},
\ee
which is very similar to, yet slightly better than, the one established in \cite{CDS} by complex variable arguments.
\end{remark}

\begin{remark}
\label{remfin}
As a consequence of Stechkin's lemma, under the assumptions of Corollary \ref{corleg}, we find that the
best $n$-term approximation polynomials
\be
u_n^L:=\sum_{\nu\in \Lambda_n^L} u_\nu L_\nu,
\ee
obtained by retaining the indices of the $n$ largest $\|u_\nu\|_V$,
satisfy the estimate 
\be
\|u-u_n^L\|_{L^2(U,V,\mu)}\lsim n^{-r},
\ee
where $r:=\frac 1 p-\frac 1 2=\frac 1 q$. There is, however, a more direct
and constructive way of retrieving this convergence rate,
namely taking instead $\Lambda_n^L$ to be the set of indices corresponding to the $n$ smallest
values of the weights $w_\nu:=a_\nu^{-1}\rho^{\nu}$ which appear in \iref{weightedell2}.
We then directly obtain that
\be
\|u-u_n^L\|_{L^2(U,V,\mu)}\leq
\sup_{\nu\notin\Lambda_n} w_\nu^{-1} \Bigl(\sum_{\nu\in\cF} w_\nu^2 \|u_\nu\|_V^2\Bigr)^{1/2} \lsim d_{n+1}^*,
\ee
where $(d_n^*)_{n\geq 1}$ is the decreasing rearrangement of the sequence $(w_\nu^{-1})_{\nu\in \cF}$.
As seen in the proof of  Corollary \ref{corleg}, this sequence 
belongs to $\ell^q(\cF)$, which implies that $d_n^*\lsim n^{-r}$ with $r:=\frac 1 q$.
\end{remark}

{\begin{remark}\label{rem:jacobi}
The results of this section can be generalized to other families of orthogonal polynomials satisfying a suitable Rodrigues' formula. For instance, when the uniform measure on $U$ is replaced by a tensor product beta measure
\be\label{betameasure}
    d\tilde\mu(y) = \bigotimes_{j\geq 1} \frac{\Gamma(\alpha_j + \beta_j + 2)}{2^{\alpha_j + \beta_j + 1}  \Gamma(\alpha_j+1) \Gamma(\beta_j+1) } (1- y_j)^{\alpha_j} (1+ y_j)^{\beta_j} dy_j,
\ee
where $\alpha_j , \beta_j > -1$ with the uniform measure as the special case $\alpha_j = \beta_j = 0$, the corresponding orthonormal polynomials are the Jacobi polynomials given by the Rodrigues' formula
\be
   P_\nu(y) = \prod_{j\geq 1} \frac{c^{\alpha_j,\beta_j}_{\nu_j} }{ \nu_j! \, 2^{\nu_j} } 
     (1 - y_j)^{-\alpha_j} ( 1+ y_j)^{-\beta_j} \,\partial^{\nu_j}_{y_j} 
       \Bigl( (y_j^2 - 1)^{\nu_j} (1-y_j)^{\alpha_j} (1+y_j)^{\beta_j} \Bigr),
\ee
where 
\be
   c^{\alpha_j, \beta_j}_{\nu_j}  := \sqrt{ \frac{ (2\nu_j + \alpha_j + \beta_j + 1)\, \nu_j!\, \Gamma(\nu_j + \alpha_j + \beta_j + 1) \,\Gamma(\alpha_j+1)\, \Gamma(\beta_j+1) }{ \Gamma(\nu_j + \alpha_j + 1)\, \Gamma(\nu_j + \beta_j + 1) \, \Gamma(\alpha_j + \beta_j + 2)} } ,
\ee
with the convention that $c^{\alpha_j, \beta_j}_{0}=1$ for any $\alpha_j , \beta_j > -1$.
Using integration by parts analogously to \eqref{rodriguesapplication}, for the Jacobi coefficients $\tilde u_\nu$ we then obtain 
\be
 \tilde u_\nu = \int_U u(y)\, P_\nu(y)\,d\tilde\mu(y) = \tilde a_\nu \int_U \frac{1}{\nu!} \partial^\nu u(y) \prod_{j\geq 1} \frac{(1-y_j^2)^{\nu_j}}{2^{\nu_j}} d\tilde\mu(y),
 \qquad \tilde a_\nu := \prod_{j \geq 1} c^{\alpha_j ,\beta_j}_{\nu_j} .
\ee
For $\rho = (\rho_j)_{j \geq 1}$ satisfying the assumptions of Theorem \ref{thm:legendreweightedsum}, proceeding exactly as in \eqref{legendreestimate} we arrive at
\be
     \sum_{\nu \in \cF} \bigl(\tilde a_\nu^{-1} \rho^\nu\|\tilde u_\nu\|_V \bigr)^2  < \infty .
\ee
Provided that the sequences $(\alpha_j)_{j\geq 1}$, $(\beta_j)_{j\geq 1}$ are such that $c^{\alpha_j,\beta_j}_{\nu_j} \leq Q(\nu_j)$ for some polynomial $Q$ independent of $j$ (which holds, for instance, when both sequences are constant), the statement of Corollary \ref{corleg} thus holds also for $\bigl(  \| \tilde u_\nu\|_V \bigr)_{\nu\in\cF}$.
\end{remark}
\begin{remark}
\label{remgen}
The same remarks as given at the end of \S 2 apply to the
generalization of Theorem \ref{thm:legendreweightedsum}, Corollary \ref{corleg}, and Remark \ref{rem:jacobi} to other types of linear PDEs with affine 
dependence in the coefficients.
\end{remark}
}

\section{Examples}\label{sec:examples}

In this section, we compare the summability properties obtained with the approach in the
present paper to those obtained with the
analysis in \cite{CDS} for various types of $(\psi_j)_{j\geq 1}$.
We show that this approach gives an improvement on Theorem \ref{firsttheo}, 
depending on the particular structure of the supports of the $\psi_j$.

\subsection{Finitely overlapping supports}\label{ssec:disj}

We first consider families $(\psi_j)_{j\geq 1}$ of functions with finitely overlapping supports, that is,
such that for any $x\in D$, 
\be
\#\{j \, : \, \psi_j(x)\neq 0\} \leq M,
\ee
for some fixed $M>0$, {where $\# S$ denotes the cardinality of the set $S$}. 
The case $M=1$ corresponds to disjoint supports, such as the family of characteristic functions 
$\psi_j=b_j\Chi_{D_j}$ with some normalizing factor $b_j$, when $(D_j)_{j\geq 1}$ is a partition of $D$. 
Another example with $M\geq 1$ is the set of Lagrange finite element basis functions
of a given order $k\geq 1$, associated to a conforming simplicial partition of $D$.

Assuming \iref{uea1}, we then define a
weight sequence $(\rho_j)_{j\geq 1}$ by
\be
\rho_j:=1+\frac {\bar a_{\min} (1-\theta)}{2M\|\psi_j\|_{L^\infty}}.
\ee
With such a choice, we find that \iref{suea} holds since, for all $x\in D$,
\be
\sum_{j\geq 1} \rho_j|\psi_j(x)| \leq \sum_{j\geq 1} |\psi_j(x)|+\frac {\bar a_{\min} (1-\theta)}{2M} \sum_{j\geq 1}\frac { |\psi_j(x)|}{\|\psi_j\|_{L^\infty}}
\leq \theta \bar a (x)+ \frac{1-\theta}{2}\bar a_{\min} \leq \delta \bar a(x),
\ee
with $\delta:=\frac {1+\theta}2<1$. As a consequence of Corollaries \ref{sumresult} and \ref{corleg} we obtain the following result.

\begin{corollary}\label{finoverlapresult}
Assume that $(\psi_j)_{j\geq 1}$ is a family of functions with finitely overlapping supports,
and that \iref{uea1} holds. If $(\| \psi_j \|_{L^\infty})_{j \geq 1} \in \ell^q(\N)$  { for some $0<q<\infty$
and if $0<p<2$ is such that $q=q(p):=\frac {2p}{2-p}$},
then $(\|t_\nu\|_V)_{\nu\in\cF}$ and $(\|u_\nu\|_V)_{\nu\in\cF}$ belong to $\ell^p(\cF)$.
In particular, best $n$-term Legendre approximations
converge in $L^2(U,V,\mu)$ with rate $n^{-s}$ where $s=\frac 1 p-\frac 1 2=\frac 1 q$.
\end{corollary}

As already observed we always have $q(p)>p$ which shows that there is in this case
a significant improvement between the summability properties of $(\|\psi_j\|_{L^\infty})_{j\geq 1}$
and those of $(\|t_\nu\|_V)_{\nu\in \cF}$ and $(\|u_\nu\|_V)_{\nu\in\cF}$, in contrast to
Theorem \ref{firsttheo}. { Note in particular that the latter
would lead to the weaker conclusion that best $n$-term Legendre approximations
converge in $L^2(U,V,\mu)$ with rate $n^{-s}$ where $s=\frac 1 q-\frac 1 2$ instead of $\frac 1 q$.}

We next give a specific example which shows that, for the Taylor coefficients,
this new result is in fact sharp. In this example, we let
$D := ]0,1[$ and $\bar a = 1$, and we consider a sequence $(D_j)_{j \geq 1}$ of disjoint intervals $D_j=]l_j,r_j[\subset D$.
Let $m_j := \frac12(l_j + r_j)$ be the midpoint of $D_j$ and $\psi_j := b_j \Chi_{[l_j, m_j]}$ with $(b_j)_{j \geq 1} \in \ell^q(\N)$, 
where $q=q(p):=\frac {2p}{2-p}$ for some $0<p<2$. We denote by
\be
h_j(x) := \max \bigl\{ 0, 1 - 2|x-m_j|/|D_j| \bigr\},
\ee
the hat function on $D_j$ centered at $m_j$ with $h_j(m_j) = 1$. We fix a sequence $(c_j)_{j \geq 1}$ such that 
$\sum_{j\geq 1}c_j^2/ |D_j| < \infty$ and choose the right hand side $f=-(\sum_{j\geq 1} c_j h_j)''\in V'$ so that
\be
t_0 = \sum_{j\geq 1} c_j h_j.
\ee
The condition on the $c_j$ ensures that $t_0 \in V$. For the particular $\nu=e_j$, the Taylor
coefficients satisfy
\be
   \int_D  t_{e_j}'  v' \,dx = - \int_D \psi_j  t_0' v'\,dx \,, \quad v \in V\,.
\ee
Testing this with $v= h_j$, by the Cauchy-Schwarz inequality we obtain
\be
   \biggl( \int_D | t_{e_j}' |^2\,dx  \biggr)^{\frac12} \biggl(  | D_j | \frac4{|D_j|^2} \biggr)^{\frac12}  \geq  | D_j | b_jc_j \frac2{|D_j|^2},
\ee
and hence 
\be
   \| t_{e_j}\|_V  \geq  \frac{ b_j c_j}{ \sqrt{| D_j |}}.
\ee
In view of the requirements $(b_j)_{j\geq 1} \in \ell^q(\N)$ and $\bigl(c_j/\sqrt{|D_j|} \bigr)_{j\geq 1} \in \ell^2(\N)$, these sequences can be chosen to ensure $(\| t_{e_j}\|_V )_{j \geq 1} \notin \ell^{\tilde p}(\N)$ for any $\tilde p< p$. 

\subsection{Arbitrary supports}\label{ssec:arbitrary}

In the general case where the supports of the $\psi_j$ are arbitrary, in particular 
for globally supported functions, the approach based on Theorem \ref{thm:weightedsum}
does not bring any specific
improvement over Theorem \ref{firsttheo} (which can be derived from it,
as observed in Remark \ref{remcor}).
One way to see this is to observe that in certain situations, 
the latter can already be sharp.
 
Consider for example the case of constant $\psi_j \equiv b_j$ with $(b_j)_{j \geq 1} \in \ell^p(\N)$ and $(b_j)_{j\geq 1} \notin \ell^{\tilde p}(\N)$ for any $\tilde p < p$,
and $\bar a=1$. For such $\psi_j$, one has
\be
u(y)=\frac 1{1+\sum_{j\geq 1} y_j b_j} u(0),
\ee
so that the Taylor coefficients are explicitly given by
\be
t_\nu=(-1)^{|\nu|} \frac {|\nu| !}{\nu!} b^{\nu} t_0.
\ee
In particular, one has  $\| t_{e_j}\|_V= b_j \| t_0\|_V$, which shows that $( \|t_\nu\|_V )_{\nu\in\cF}\notin \ell^{\tilde p}(\cF)$ for any $\tilde p < p$. 
A similar, yet more technical, computation shows that the same holds for the Legendre coefficients.
Therefore, in this case of completely overlapping supports of the $\psi_j$, the proposed new bounds cannot give an improvement
over Theorem \ref{firsttheo}.

\subsection{Wavelets}\label{ssec:wavelets}

Let us now turn to the case of diffusion coefficients parametrized by a wavelet basis, that is,
\be
 a(y) = \bar a + \sum_{\lambda} y_\lambda \psi_\lambda,
\ee
where 
\be
\| \psi_\lambda \|_{L^\infty} = c_l:=C 2^{- \alpha l},\quad  |\lambda|=l,
\label{dec}
\ee
for some $\alpha >0$,
as discussed in the introduction. Note that, when ordering the wavelet basis from coarse to fine scales, the 
resulting system $(\psi_j)_{j\geq 1}$ has then the algebraic behaviour
\be
\|\psi_j\|_{L^\infty} \sim j^{-\alpha/d}.
\ee
Assuming \iref{uea1}, and for an arbitrary $0<\beta<\alpha$ we define a weight sequence $(\rho_\lambda)$ by
\be
\rho_{\lambda}:=1+   \frac{\bar a_{\min}  (1-\theta) }{2CM(1-2^{\beta-\alpha}) } 2^{\beta |\lambda|},
\ee
where $M$ and $C$ are the constants in \iref{finiteoverlap} and \iref{dec}. With such a choice, we find that
\iref{suea} holds since, for all $x\in D$,
\be
\sum_{\lambda}\rho_{\lambda}|\psi_\lambda(x)|
\leq \theta \bar a(x)+\sum_{\lambda} \frac{\bar a_{\min}  (1-\theta) }{2CM(1-2^{\beta-\alpha}) } 2^{\beta |\lambda|}|\psi_j(x)|
\leq \theta \bar a(x)+\frac{\bar a_{\min}  (1-\theta) }{2} \leq \delta \bar a(x),
\ee
where $\delta:=\frac {1+\theta}2<1$. After the same reordering as for the wavelet basis, we find that
\be
\rho_j \sim  j^{\beta/d}.
\ee
Therefore, as a consequence of Corollaries \ref{sumresult} and \ref{corleg} we obtain the following result.

\begin{corollary}\label{waveletresult}
Assume that $(\psi_j)_{j\geq 1}$ is a wavelet basis with normalization \iref{dec}
and that \iref{uea1} holds. If $(\| \psi_j \|_{L^\infty})_{j \geq 1} \in \ell^q(\N)$ for some $q<q(p):= \frac{2p}{2-p}$,
then $(\|t_\nu\|_V )_{\nu\in \cF}$ and $(\|u_\nu\|_V )_{\nu\in \cF}$ belong to $\ell^p(\cF)$.
In particular, best $n$-term Legendre approximations
converge in $L^2(U,V,\mu)$ with rate $n^{-s}$ for all $s<\frac 1 q$.
\end{corollary}

As already mentioned, if we use sufficiently smooth wavelets, the decay property \iref{dec} is equivalent
to the property that $a(y)$ is in the Besov space $B^\alpha_\infty(L^\infty(D))$,
which for non-integer $\alpha$ coincides with the H\"older space $C^\alpha(D)$, for all $y\in U$. Thus,
we also infer from Corollary \ref{corleg} that if this holds
for some $\alpha>0$, best $n$-term Legendre approximations
converge in $L^2(U,V,\mu)$ with rate $n^{-s}$ for all $s<\alpha/d$.

\section{Numerical illustrations}

In the following numerical tests, we consider three different cases of parametrized diffusion problems on $D=]0,1[$. In each of these cases, the parameter dependence is expressed in terms of a different representative type of function system $( \psi_j)_{j\geq 1}$ that corresponds to one of the three types of examples considered in Section \ref{sec:examples}.

We give results both for Taylor and for Legendre coefficients. The Taylor coefficients are computed using the alternating greedy Taylor algorithm described in \cite[\S 7.2]{CD}. For the Legendre coefficients we use an adaptive stochastic Galerkin method using the approach described in \cite{G}, where multi-index sets are refined based on approximate evaluation of the residual of the continuous problem. Both methods are guaranteed to converge, but in the Taylor case one needs to resort to a heuristic stopping criterion. For each of the tests, a suitably adjusted single fixed finite element discretization using P2 elements is used for the spatial dependence.

In all of the following examples, we take $\bar{a} := 1$. We also fix $\theta \in ]0,1[$ and choose the functions $\psi_j$ in the examples such that this $\theta$ satisfies \eqref{thetadef}.

In the following discussion of the methodology for evaluating the results, we always refer only to the Taylor coefficients $(t_\nu)_{\nu\in\cF}$, with the understanding that the same considerations apply to the Legendre coefficients $(u_\nu)_{\nu\in\cF}$.
In each of the tests, the aim is to numerically estimate the limiting summability exponent
\be\label{barp}
 \bar p := \inf \bigl\{ \,p > 0 \:\colon\:  (\|t_\nu\|_V)_{\nu\in\cF}  \in \ell^p(\cF) \bigr\}.
\ee
To this end, we introduce the decreasing rearrangement $(t^*_n)_{n\geq 1}$ of $(\|t_\nu\|_V)_{\nu\in\cF}$. Then $(\|t_\nu\|_V)_{\nu\in\cF} \in \ell^p(\cF)$ implies that for some $C>0$, one has $t^*_n \leq C n^{-1/p}$,
and conversely, if $t^*_n \leq C n^{-1/q}$ for some $C,q>0$, then  $(\|t_\nu\|_V)_{\nu\in\cF} \in \ell^p(\cF)$ for any $p>q$. As a consequence, 
\be\label{pinf}
  \bar p = \inf \bigl\{ \, p>0 \: \colon \:  \sup_{n\in \N} n^{\frac1p} t^*_n  < \infty \bigr\},
\ee 
or in other words, $\bar p$ can be determined from the asymptotic decay rate of the values $t^*_n$.
As estimates for the largest $s>0$ such that $\sup_n n^s t^*_n $ is finite, we consider the values 
\be\label{si}
       s_i: = \log_2 (t^*_{2^{i-1}}) - \log_2 (t^*_{2^{i}}),
\ee
for $i=1,2,\dots$.
In view of \eqref{pinf}, if the sequence $(t^*_n )_{n\geq 1}$ decays asymptotically at an algebraic rate, for sufficiently large $i$ we can thus expect $s_i$ to approximate $\bar p^{-1}$.

\subsection{Parametrization by disjoint inclusions}\label{secnumdisjoint}

In the first test, we choose a family $\{ D_j \}_{j\geq 1}$ of disjoint open intervals in $D$ and a $\beta >0$, and define
\be
  \psi_j := \theta j^{-\beta} \Chi_{D_j}.
\ee
Note that although this does not enter into any of the decay estimates available at this point, the concrete example in \S \ref{ssec:disj} suggests that the decay of the inclusion sizes $|D_j|$ has an impact on the summability of $(\|t_\nu\|_V)_{\nu\in\cF}$ and $(\|u_\nu\|_V)_{\nu\in\cF}$. Indeed, it can also observed numerically that faster decay of $|D_j|$ leads to improved summability. To remove this effect in the tests, we therefore choose $D_j$ such that the decay of $|D_j|$ is as slow as possible
while still allowing $D_j$ to partition $D$. To this end, we define 
\be
x_0 := 0\quad {\and} \quad x_j := c \sum_{k=1}^j k^{-1} \log^{-2}(1+k), \quad j\geq 1,
\ee 
with $c$ such that $\lim_{j\to \infty} x_j = 1$, and set $D_j = ] x_{j-1}, x_j [$. Since $(\| \psi_j \|_{L^\infty})_{j \geq 1} \in \ell^q(\N)$ for all $q > \frac1\beta$, by Corollary \ref{finoverlapresult} we expect that $(\|t_\nu\|_V)_{\nu\in\cF}$ and $(\|u_\nu\|_V)_{\nu\in\cF}$ belong to $\ell^p(\cF)$ for any $p > (\beta + \frac12)^{-1}$.

\begin{figure}[tp]\centering
\begin{tabular}{cc}
\includegraphics[width=7cm]{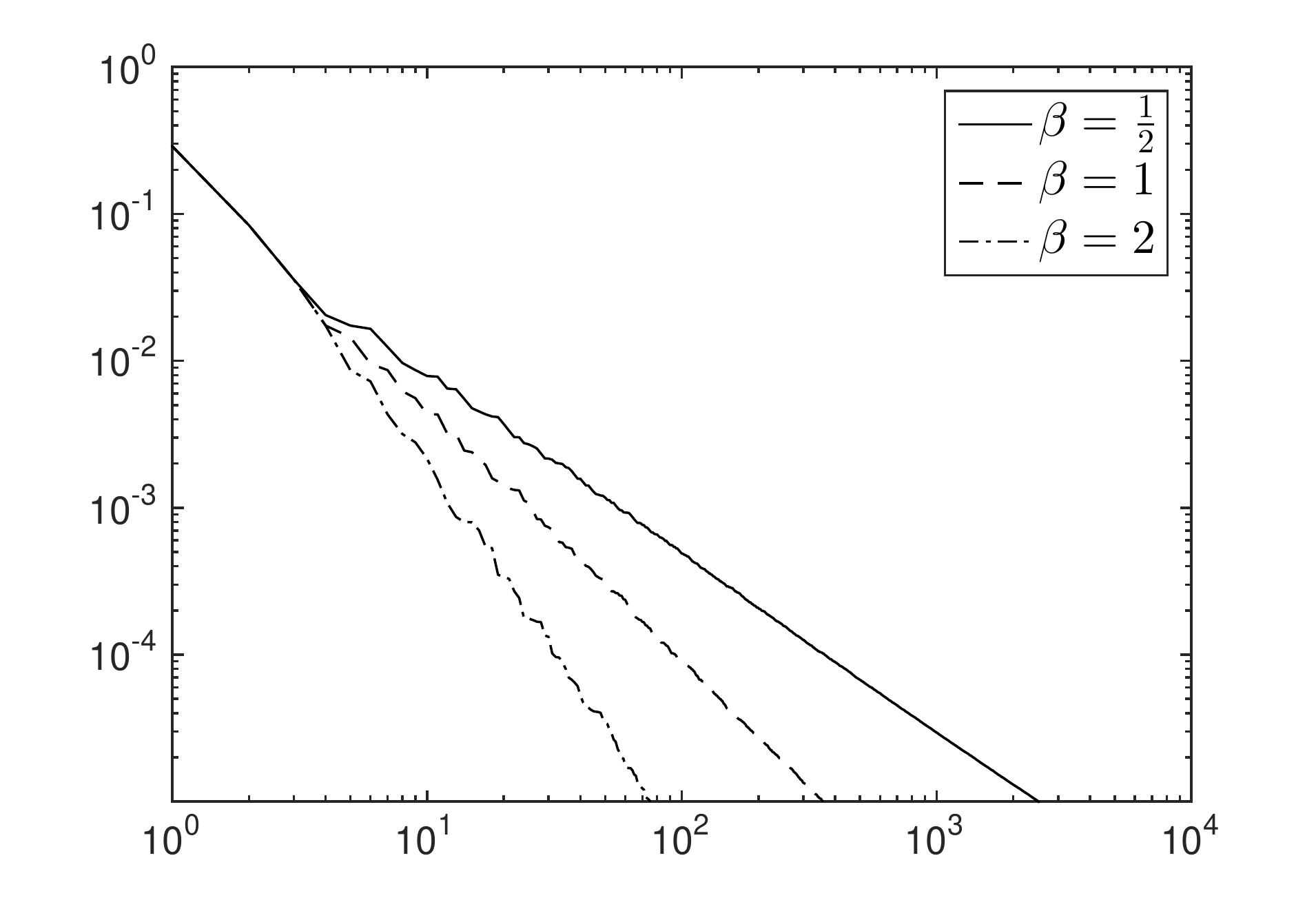}  & \includegraphics[width=7cm]{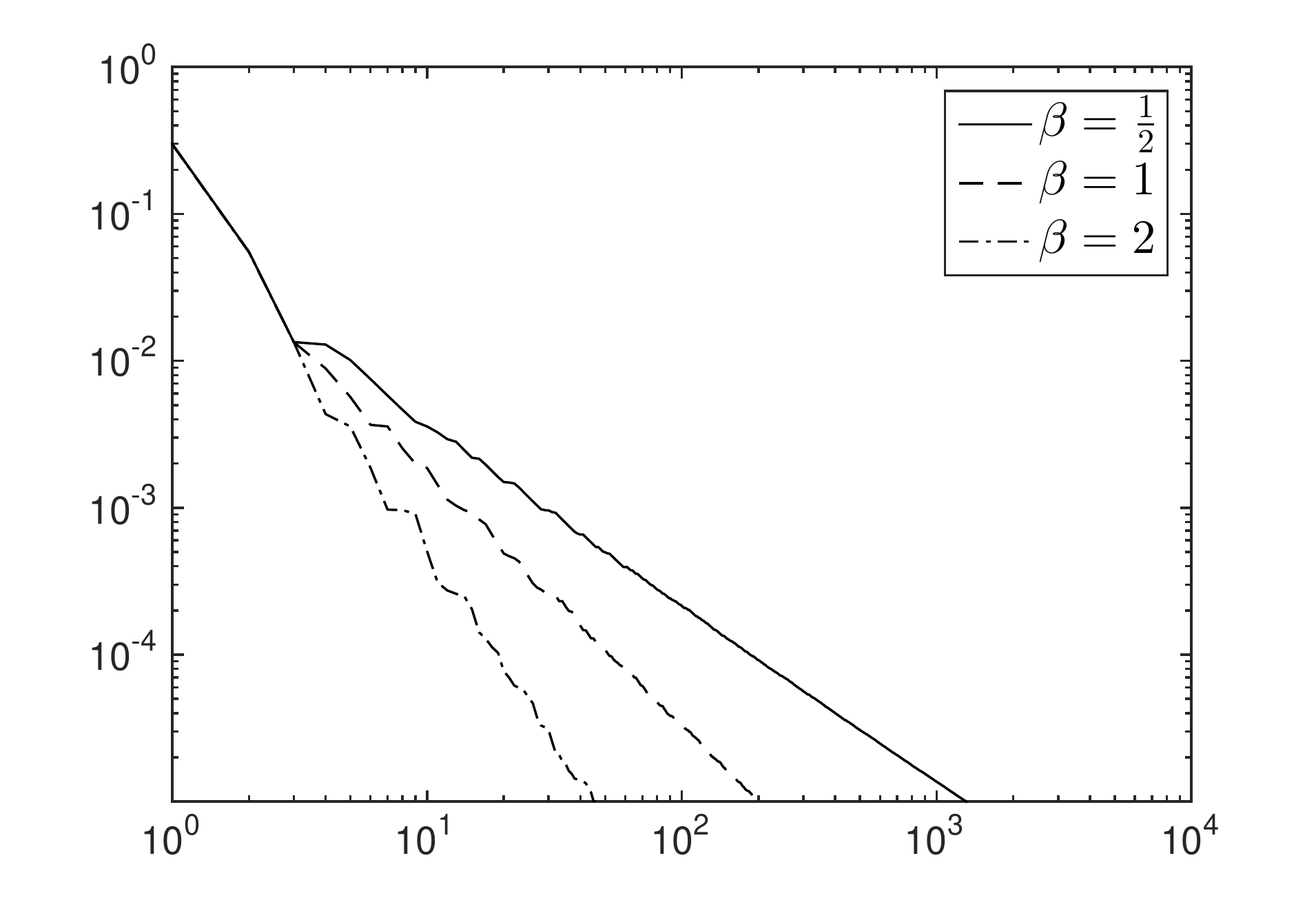} \\[-6pt]
\small Taylor   & \small Legendre
\end{tabular}
\caption{Disjoint inclusions: ordered norms of Taylor coefficients $t_\nu$ and Legendre coefficients $u_\nu$, for $\theta=\frac12$ and $\beta = \frac{1}{2}, \,1,\, 2$.}
\label{disjointfig}
\end{figure}

\begin{table}[tp]\centering\small
\begin{tabular}{l|rrr||rrr}
  & \multicolumn{3}{c}{Taylor}  & \multicolumn{3}{c}{Legendre} \\ \hline
      &   $\beta=2$     &     $\beta=1$     &     $\beta=\frac{1}{2}$    &      $\beta=2$     &     $\beta=1$     &     $\beta=\frac{1}{2}$\\ \hline
$s_6$    &  2.563   &   1.730   &   1.225   &   2.476   &   1.789  &   1.302  \\
$s_7$    &  2.708   &   1.731   &   1.274   &   2.578   &   1.786  &   1.235  \\
$s_8$    &  2.481   &   1.726   &   1.211   &   2.601   &   1.701  &   1.212  \\
$s_9$    &  2.574   &   1.706   &   1.235   &   2.514   &   1.661  &   1.200  \\
$s_{10}$ &  2.439   &   1.650   &   1.196   &   2.543   &   1.660  &   1.169 \\
$s_{11}$ &  2.477   &   1.643   &   1.175   &   2.507   &   1.642  &   1.160 \\ \hline
${\bar p}^{-1}$  & 2.500 & 1.500  & 1.000  & 2.500 & 1.500 & 1.000
\end{tabular}
\caption{Disjoint inclusions: decay rates of coefficient norms, with $s_i$ as in \eqref{si}, compared to limiting value ${\bar p}^{-1}=\beta+\frac12$ expected by Corollary \ref{finoverlapresult}. }
\label{disjointtab}
\end{table}
The values of the decreasing rearrangements of these sequences for $\beta = \frac12, 1, 2$, where in each case $\theta=\frac12$, are compared in Figure \ref{disjointfig}. In Table \ref{disjointtab}, the empirically determined decay rates are compared to the theoretical prediction for $\bar p^{-1}$. We observe almost the same decay behavior for Taylor and Legendre coefficients, and in each case the empirical rates indeed approach $\bar p^{-1}$.

\subsection{Parametrization by a Fourier expansion}\label{secnumfourier}

We next consider a parametrization by the globally supported Fourier basis
\be\label{sineexample}
   \psi_j (x):= \theta c j^{-\beta} \sin (j \pi x),
\ee
for some $\beta>1$, with the normalization constant $c := \bigl(\sum_{j\geq 1} j^{-\beta} \bigr)^{-1}$. We thus have $(\| \psi_j \|_{L^\infty})_{j\geq 1} \in \ell^p(\N)$ for all $p > \frac1\beta$.
In view of the discussion in \S \ref{ssec:arbitrary}, due to Theorem \ref{firsttheo} we expect that $(\|t_\nu\|_V)_{\nu\in\cF}$ and $(\|u_\nu\|_V)_{\nu\in\cF}$
belong to $\ell^p(\cF)$ for such $p$. 

\begin{figure}[tp]\centering
\begin{tabular}{cc}
\includegraphics[width=7cm]{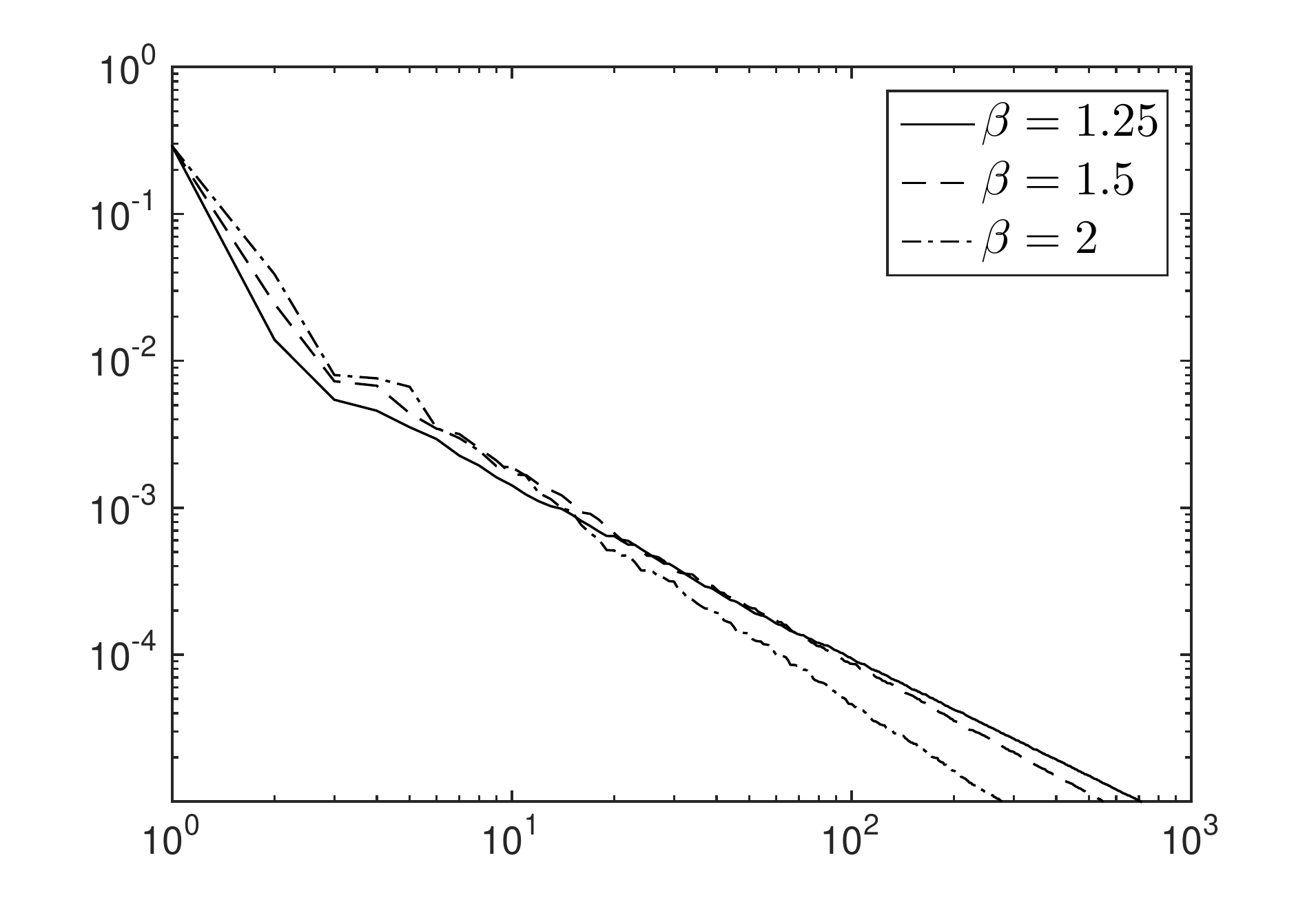}  & \includegraphics[width=7cm]{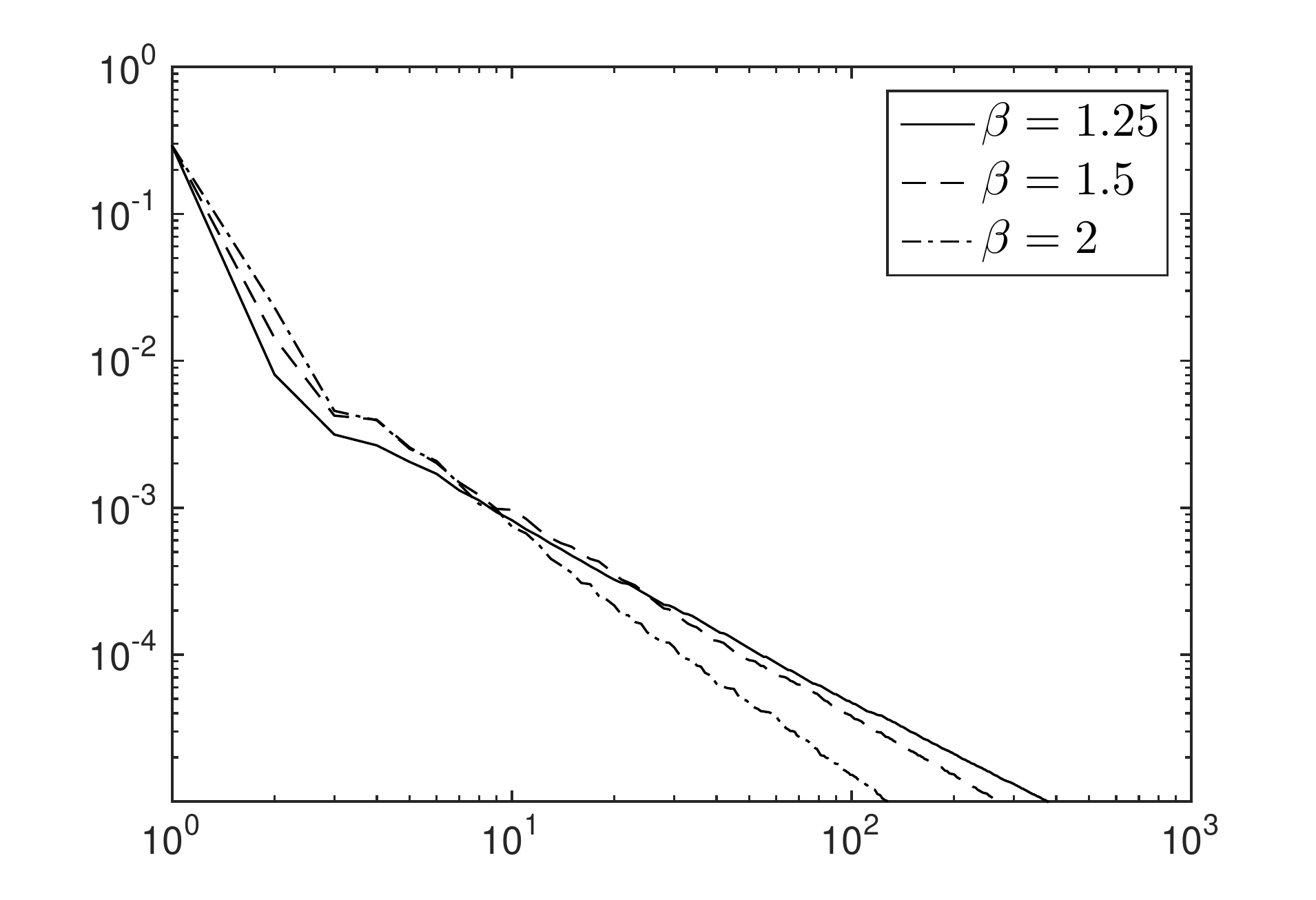} \\[-6pt]
\small Taylor   & \small Legendre
\end{tabular}
\caption{Fourier expansion: ordered norms of Taylor coefficients $t_\nu$ and Legendre coefficients $u_\nu$, for $\theta=\frac12$ and $\beta = 1.25, 1.5, 2$.}
\label{fourierfig}
\end{figure}
\begin{table}[tp]\centering\small
\begin{tabular}{l|rrr||rrr}
  & \multicolumn{3}{c}{Taylor}  & \multicolumn{3}{c}{Legendre} \\ \hline
         & $\beta=2$  &  $\beta=1.5$  &     $\beta=1.25$  &    $\beta=2$    &  $\beta=1.5$  &    $\beta=1.25$  \\ \hline
$s_6$    &    1.452    &      1.165    &     1.250        &     1.593       &     1.294    &       1.250 \\
$s_7$    &    1.619    &      1.320    &     1.092        &     1.682       &     1.353    &       1.154 \\
$s_8$    &    1.495    &      1.278    &     1.147        &     1.597       &     1.337    &       1.192 \\
$s_9$    &    1.515    &      1.257    &     1.141        &     1.632       &     1.338    &       1.187  \\
$s_{10}$ &    1.533    &      1.270    &     1.143        &     1.637       &     1.341    &       1.173 \\
$s_{11}$ &    1.515    &      1.258    &     1.143        &     1.639       &     1.327    &       1.191  \\ \hline
${\bar p}^{-1}$  & 2.000 & 1.500  & 1.250  & 2.000 & 1.500 & 1.250
\end{tabular}
\caption{Fourier expansion: decay rates of coefficient norms, with $s_i$ as in \eqref{si}, compared to limiting value $\bar p^{-1}=\beta$ expected by Theorem \ref{firsttheo}. }
\label{fouriertab}
\end{table}
The results for $\theta = \frac12$ and $\beta = 1.25, 1.5, 2$ are shown in Figure \ref{fourierfig} and Table \ref{fouriertab}. Here we observe that especially for larger values of $\beta$, the empirically observed rates $s_i$ do not come very close to the theoretically guaranteed limiting value ${\bar p}^{-1}$ within the considered range of coefficients. This indicates that the asymptotic behavior emerges only very late in the expansions. 

Note that this observation is consistent with the results obtained in \cite{G}, where the above example \eqref{sineexample} with $\beta=2$ and $\theta=\frac12$ is considered as a numerical test. There, a decay rate close to $1$ is observed for the $L^2$ error of a Legendre expansion (with fixed spatial grid), corresponding to a decay rate of the coefficient norms close to $1.5$ as obtained here.

\begin{figure}[tp]
\hspace{1.5cm}\begin{minipage}[b]{7cm}
\includegraphics[width=7cm]{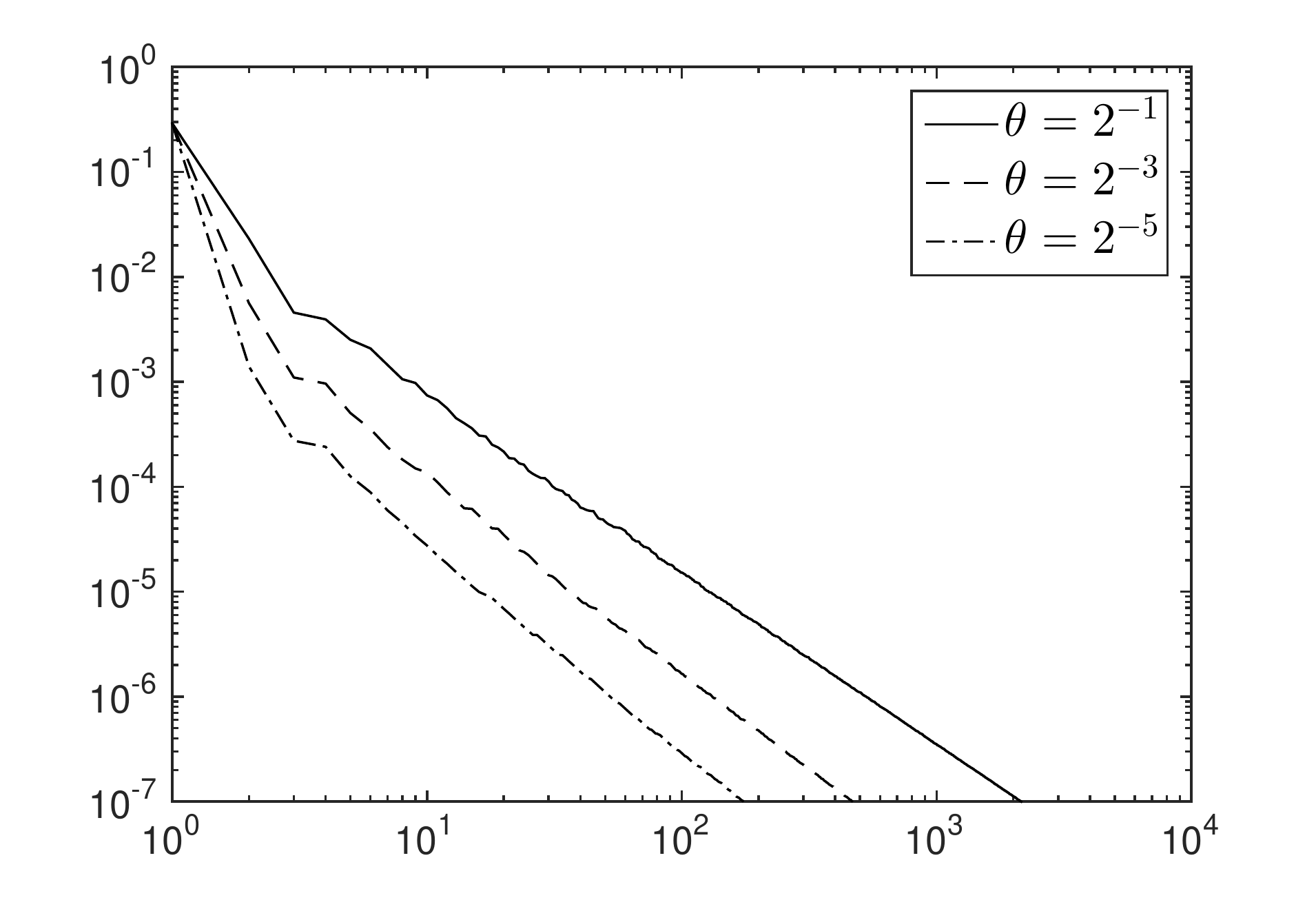}
\end{minipage}\raisebox{2.1cm}{\begin{minipage}[b]{8cm}\small
\begin{tabular}{l|rrr}
         &    $\theta=2^{-1}$ &  $\theta=2^{-3}$   &  $\theta=2^{-5}$  \\ \hline
$s_6$    &    1.593              &  1.876                 &   2.000  \\
$s_7$    &    1.682              &  1.767                 &   1.872  \\
$s_8$    &    1.597              &  1.822                 &   1.908  \\
$s_9$    &    1.632              &  1.813                 &   1.905  \\
$s_{10}$ &    1.637              &  1.813                 &   1.898  \\
$s_{11}$ &    1.639              &  1.814                 &   1.921   \\ \hline
${\bar p}^{-1}$ & 2.000  & 2.000 & 2.000
\end{tabular}
\end{minipage} }
\caption{Fourier expansion: ordered norms of Legendre coefficients $u_\nu$, for $\beta=2$ and $\theta = \frac{1}{2}, \frac{1}{2^3},\frac1{2^5}$, and corresponding decay rates.}
\label{fouriertheta}
\end{figure}

It turns out that the observed decay rates are in fact also influenced by the value of $\theta$: as shown in Figure \ref{fouriertheta}, for smaller values of $\theta$, the $s_i$ are closer to the limiting value already within the considered range.

\subsection{Parametrization by a Haar wavelet expansion}

As a final example, we return to the wavelet parametrization of $a$ with a levelwise decay \iref{dec}.
Here we use the Haar wavelet, generated from $h := \Chi_{[0,\frac12[} - \Chi_{[\frac12,1[}$, such that
\be
\psi_\lambda(x)=c_l h(2^lx-k), \quad \lambda=(l,k), \quad l\geq 0, \quad k= 0, \ldots, 2^l-1,
\ee
and we set
\be
c_l := \theta\, (1 - 2^{-\alpha}) \,2^{-\alpha l}
\ee
for a fixed $\alpha > 0$. Since, after reordering, we have $\|\psi_j\|_{L^\infty}\sim j^{-\alpha}$ and therefore $(\| \psi_j \|_{L^\infty})_{j \geq 1} \in \ell^q(\N)$ for all $q > \frac1\alpha$, by Corollary \ref{waveletresult} we expect that $(\|t_\nu\|_V)_{\nu\in\cF}$ and $(\|u_\nu\|_V)_{\nu\in\cF}$ belong to $\ell^p(\cF)$ for any $p > (\alpha + \frac12)^{-1}$.

\begin{figure}[tp]\centering
\begin{tabular}{cc}
\includegraphics[width=7cm]{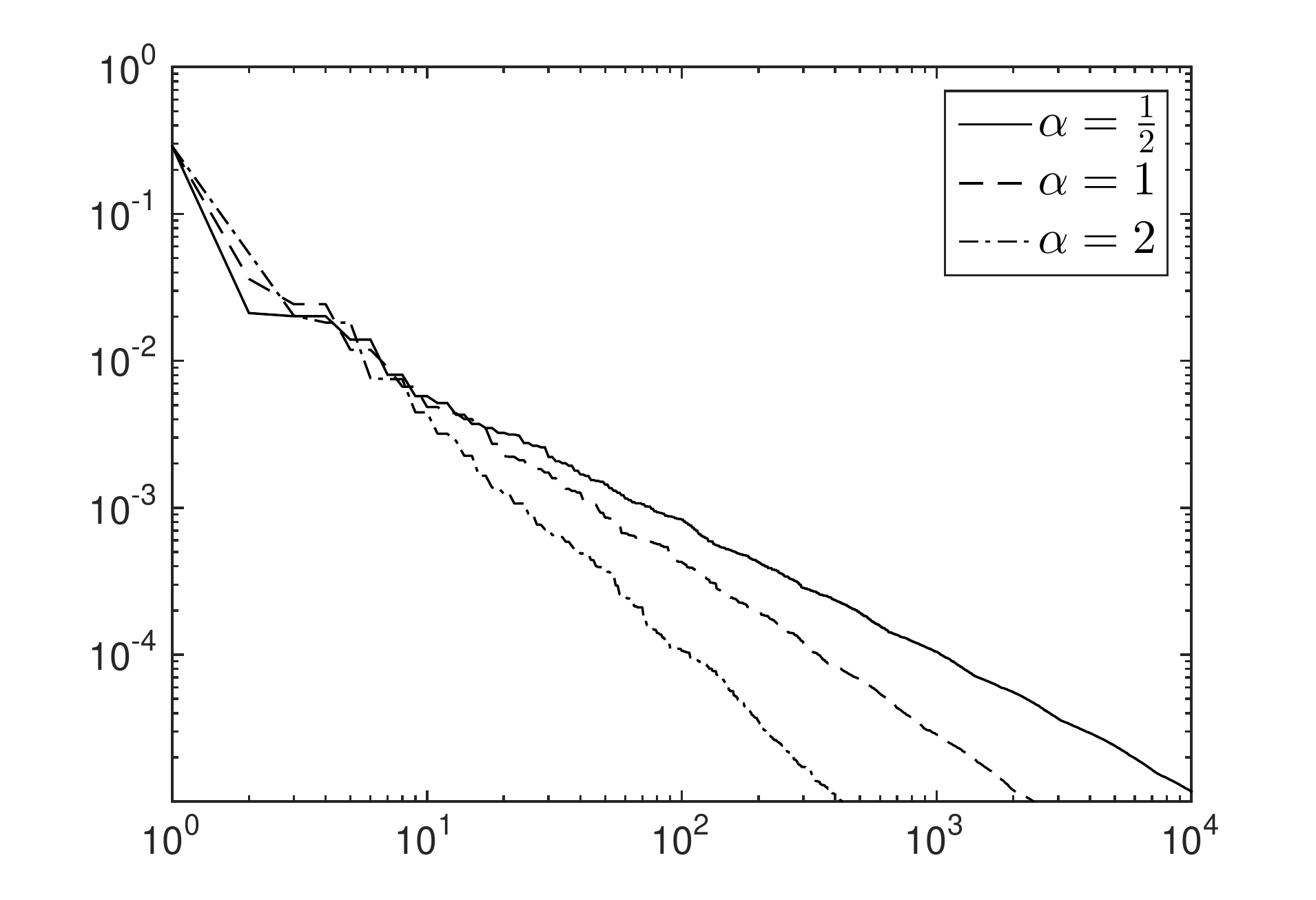}  & \includegraphics[width=7cm]{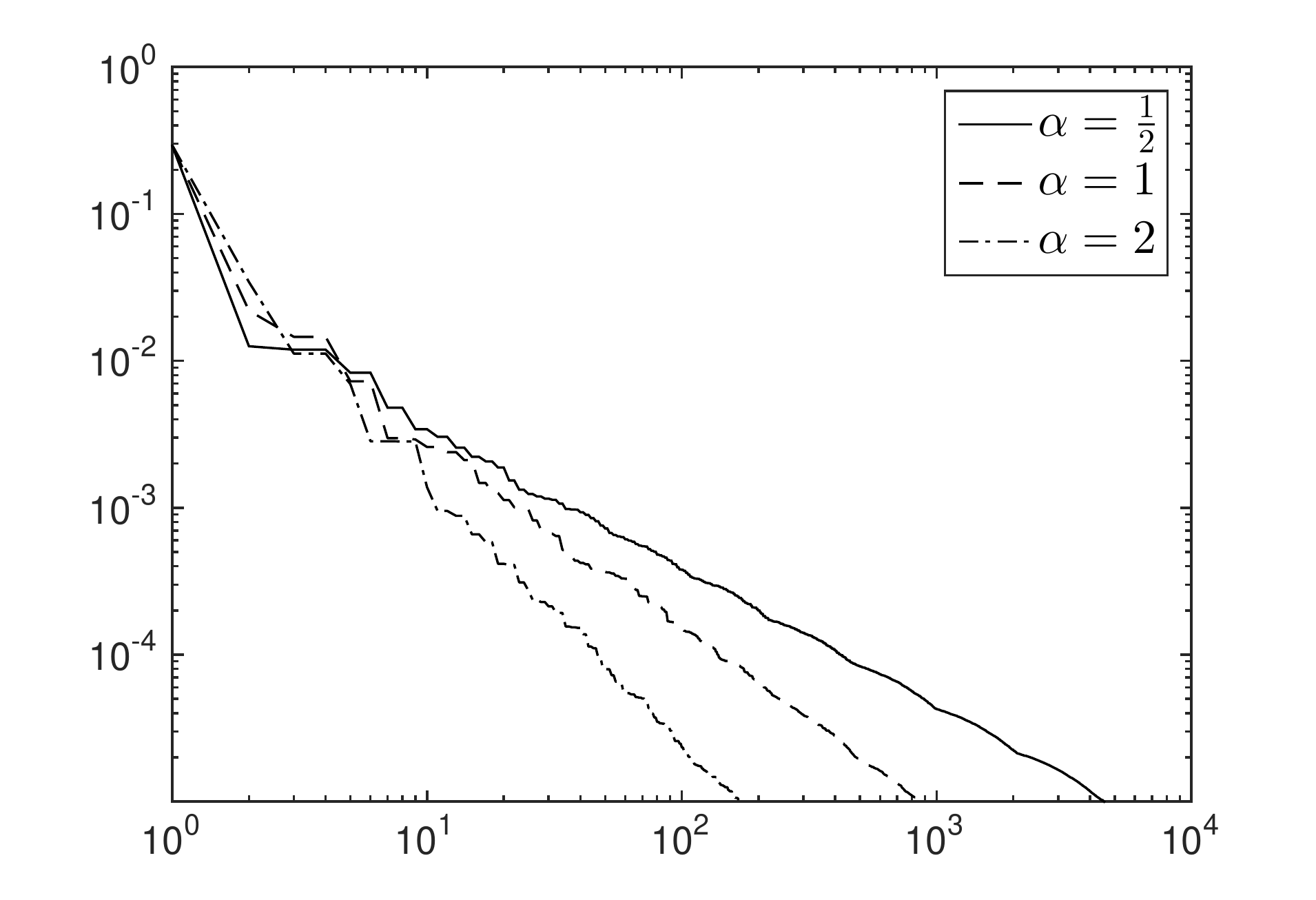} \\[-6pt]
\small Taylor   & \small Legendre
\end{tabular}
\caption{Haar expansion: ordered norms of Taylor coefficients $t_\nu$ and Legendre coefficients $u_\nu$, for $\theta=\frac12$ and $\alpha = \frac{1}{2}, \,1,\, 2$.}
\label{haarfig}
\end{figure}

\begin{figure}[tp]
\hspace{1.5cm}\begin{minipage}[b]{7cm}
\includegraphics[width=7cm]{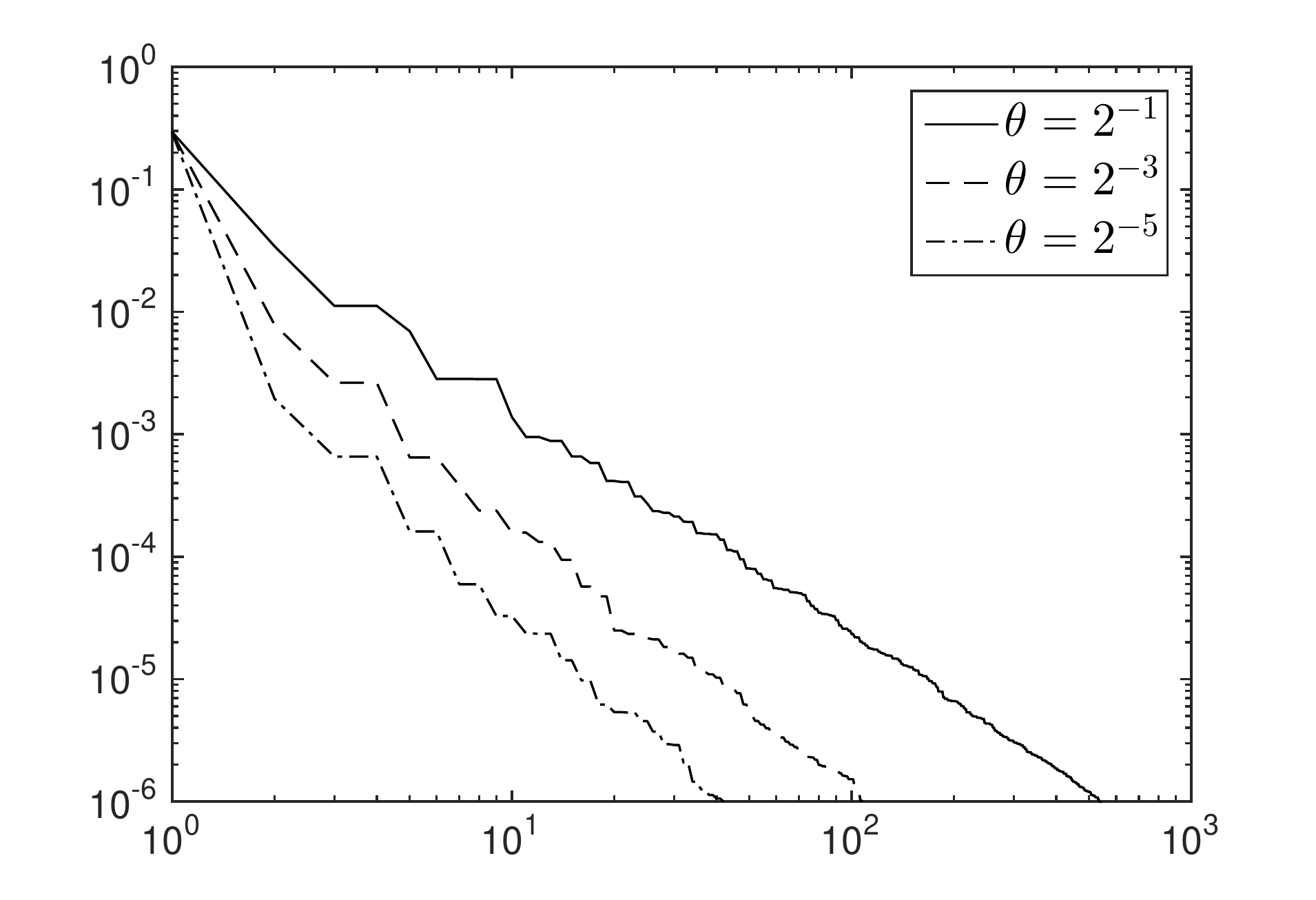}
\end{minipage}\raisebox{2.1cm}{\begin{minipage}[b]{8cm}\small
\begin{tabular}{l|rrr}
         & $\theta=2^{-1}$  & $\theta=2^{-3}$ & $\theta=2^{-5}$  \\ \hline
$s_6$    & 1.853            &  2.388          &   2.123  \\
$s_7$    & 1.779            &  2.082          &   2.175  \\
$s_8$    & 1.874            &  2.226          &    2.347  \\
$s_9$    & 1.913            &  2.056          &   2.410  \\
$s_{10}$ & 1.909            &  2.238          &   2.321  \\
$s_{11}$ & 2.037            &  2.196          &   2.396   \\ \hline
$\bar p^{-1}$ & 2.500  & 2.500 & 2.500
\end{tabular}
\end{minipage} }
\caption{Haar expansion: ordered norms of Legendre coefficients $u_\nu$, for $\alpha=2$ and $\theta = \frac{1}{2}, \frac{1}{2^3},\frac1{2^5}$, and corresponding decay rates.}
\label{haartheta}
\end{figure}

The results for $\theta = \frac12$ and $\alpha = \frac12, 1, 2$ are given in Figure \ref{haarfig} and Table \ref{haartab}. We again observe very similar decay for Taylor and Legendre expansions. 
Similarly to the observations made in Section \ref{secnumfourier}, for larger $\alpha$, the empirical rates $s_i$ do not come very close to the expected asymptotic limit $\bar p^{-1}$ within the considered range of coefficients. As shown in Figure \ref{haartheta}, the decay rates again approach $\bar p^{-1}$ more quickly for smaller values of $\theta$.
\begin{table}[tp]\centering\small
\begin{tabular}{l|rrr||rrr}
  & \multicolumn{3}{c}{Taylor}  & \multicolumn{3}{c}{Legendre} \\ \hline
         & $\alpha=2$   & $\alpha=1$     &   $\alpha=\frac12$ &  $\alpha=2$ & $\alpha=1$  & $\alpha=\frac12$ \\ \hline
$s_6$    &   1.450      &      1.301     &  0.927             &  1.853      &  1.165      &   0.947  \\
$s_7$    &   1.569      &     0.993      &  0.878             &  1.779      &  1.339      &   0.939  \\
$s_8$    &   1.794      &      1.122     &  0.803             &  1.874      &  1.275      &   0.953  \\
$s_9$    &   1.633      &      1.186     &  0.866             &  1.913      &  1.330      &   0.949  \\
$s_{10}$ &   1.799      &      1.225     &  0.872             &  1.909      &  1.247      &   0.961  \\
$s_{11}$ &   1.866      &      1.266     &  0.903             &  2.037      &  1.268      &   0.958 \\ \hline
$\bar p^{-1}$  & 2.500 & 1.500  & 1.000  & 2.500 & 1.500 & 1.000
\end{tabular}
\caption{Haar expansion: decay rates of coefficient norms, with $s_i$ as in \eqref{si}, compared to limiting value $\bar p^{-1}=\alpha+\frac12$ expected by Corollary \ref{waveletresult}. }
\label{haartab}
\end{table}

In summary, these numerical results support the conjecture that the summability estimate in Corollary \ref{waveletresult} for wavelet expansions is in fact sharp, similarly to Corollary \ref{finoverlapresult} for disjoint inclusions, whose sharpness we have established by an analytical example.

\end{document}